\documentclass[11pt]{article}
\usepackage{amsfonts}
\usepackage{bbm}
\usepackage{color}
\usepackage[toc,page,title,titletoc,header]{appendix}
\usepackage{amssymb, graphicx, amsmath}
\usepackage{amsthm,cite,indentfirst}
\usepackage{subfigure}
\usepackage{graphicx}
\usepackage{url}
\usepackage{multirow}
\usepackage{float}
\usepackage{hyperref}
\usepackage{helvet}
\usepackage{courier}
\usepackage{graphics,epsfig,graphicx,subfigure}
\usepackage{epstopdf}
\usepackage[table]{xcolor}
\usepackage{helvet}
\usepackage[table]{xcolor}
\usepackage{courier}
\usepackage{booktabs}
\usepackage{multirow}
\usepackage{dsfont}
\usepackage{pgf}
\usepackage{tikz}
\usepackage[ruled]{algorithm2e}
\usepackage[justification=centering]{caption}
\usepackage{threeparttable}
\usepackage{rotating}
\usepackage{wrapfig}
\usepackage{mathrsfs}
\allowdisplaybreaks[4]

\usetikzlibrary{arrows,decorations.pathmorphing,backgrounds,positioning, fit, petri, automata}

\tikzset{elegant/.style={smooth,thick,samples=50,line width=1.2pt}}
\tikzset{eaxis/.style={->,>=stealth}}

\renewcommand\arraystretch{1.2}

\newtheorem{theorem}{Theorem}[section]

\newtheorem{lemma}{Lemma}[section]
\newtheorem{proposition}{Proposition}[section]

\newtheorem{ex}{Example}[section]
\newtheorem{remark}{Remark}[section]

\newtheorem{assumption}{Assumption}[section]
\newtheorem{cor}{Corollary}[section]

\title{Variance-based Bregman extragradient algorithm with line search for solving stochastic variational inequalities\thanks{This work was supported by the NSF of Chongqing (cstc2021jcyj-msxmX0721, cstc2018jcyjAX0119), the Education Committee Project Research Foundation of Chongqing (KJZDK201900801) and the Innovation Project for Graduate Students of Chongqing (CYS22629).}}

\date{\today}

\author{Xian-Jun Long\thanks{\baselineskip 9pt College of Mathematics and Statistics, Chongqing Technology and Business University, Chongqing 400067, P.R.China. Corresponding author. Email: xianjunlong@ctbu.edu.cn.},
\ Yue-Hong He\thanks{\baselineskip 9pt College of Mathematics and Statistics, Chongqing Technology and Business University, Chongqing 400067, P.R.China. Email: heyuehong1111@163.com.}
\ and  \ Nan-Jing Huang\thanks{\baselineskip 9pt Department of Mathematics, Sichuan University, Chengdu, Sichuan 610064, P.R.China. Email: nanjinghuang@hotmail.com.}}

\begin{document}
\maketitle

\begin{abstract}
\noindent The main purpose of this paper is to propose a variance-based Bregman extragradient algorithm with line search for solving stochastic variational inequalities, which is robust with respect an unknown Lipschitz constant.
We prove the almost sure convergence of the algorithm by a more concise and effective method instead of using the supermartingale convergence theorem. Furthermore, we obtain not only the convergence rate $\mathcal{O}(1/k)$ with the gap function when $X$ is bounded, but also the same convergence rate in terms of the natural residual function when $X$ is unbounded. Under the Minty variational inequality condition, we derive the iteration complexity $\mathcal{O}(1/\varepsilon)$ and the oracle complexity $\mathcal{O}(1/\varepsilon^2)$ in both cases. Finally, some numerical results demonstrate the superiority of the proposed algorithm.\\
\end{abstract}


\baselineskip=16pt
\parskip=2pt

\section{Introduction}\label{sec:Introduction}
\setcounter{equation}{0}

\noindent Consider a probability space $(\Omega,\mathcal{F},\mathbb{P})$, denote by $\mathbb{E}$ the expectation with respect to the probability measure $\mathbb{P}$. Let $f:\mathbb{R}^n\times \Xi\to \mathbb{R}^n$ be a measurable function  with respect to a random variable $\xi:\Omega\to \Xi$, whose distribution is $\bf{P}$ so that $\textbf{P}(A):=\mathbb{P}(\omega\in A)$ for any $A\in \mathcal{F}$. Let $X\subseteq\mathbb{R}^n$ be a closed convex set. The stochastic variational inequality problem (SVI) is to find a vector $x^*\in X$ such that
\begin{equation}\label{equ:1:1}
\mathbb{E}[f(x,\xi(\omega))]^T(x-x^*)\geq0,\; \forall x\in X,
\end{equation}
which has many applications including in economics, game theory, modern science and engineering fields; see, e.g., [1-5]. In what follows, we use $\xi$ to denote $\xi(\omega)$ and $F(x):=\mathbb{E}[f(x,\xi(\omega))]=\int_{\Omega}f(x,\xi(\omega))d\textbf{P}(\omega)$. We also denote by $X^*$ the set of solutions of the SVI.

Note that if the operator $F$ defined in \eqref{equ:1:1} can be accurately calculated, then the SVI reduces to the deterministic variational inequality (VI), i.e., finding $x\in X$ such that $F(x)^T(y-x)\geq 0,\;\forall y\in X$, which is fundamental in a broad range of mathematical and applied sciences [6-7].
However, the SVI is very difficult to be transformed into the VI in general either because of expensive computation of the involved integration or unknown probability distribution of $\xi$, as well as no closed form of $f(\cdot,\xi)$. Therefore, the existing numerical methods for the VI may not be applied to solve the SVI directly. Recently, some effective approaches to solve the SVI have been proposed, for example: sample average approximation (SAA) approach and stochastic approximation (SA) approach. In this paper, we mainly focus on the SA approach for its significant efficiency in dealing with large-scale problems in the fields of statistics, machine learning and so on.

\subsection{Extragradient schemes for the SVI}\label{sec:1.1}
\noindent The SA approach was initially introduced by Robbins and Monro \cite{RM} for solving the stochastic root-finding problem, whose strategy is to replace the exact mean operator $F$ along the iteration by a random sample of $f$. Here, the samples are accessed in an interior and online fashion that generates the stochastic error $f(x,\xi)-F(x)$. More specifically, the first SA-based projection method for the SVI was given by Jiang and Xu \cite{JX} and they proved the convergence when $F$ is strongly monotone and Lipschitz continuous. Koshal et al. \cite{KNS1} proposed a stochastic iterative Tikhonov regularization method for the SVI, where the mapping is  monotone and Lipschitz continuous. Some related results can be found in \cite{CS,YNS,CS1}.

Recently, the SA approach was extended from single sample to multisamples for pursuing better convergence rate and complexity results, which has attracted a lot of attention from the research community. Particularly, one of the most popular method for the SVI is the following stochastic extragradient algorithm introduced by Iusem et al. \cite{I2017}:
\begin{equation}\label{equ:1:2}
\left\{
\begin{array}{ll}
x_{k+1/2}:=\Pi_X(x_k-\gamma \frac{1}{N_k}\sum^{N_k}_{j=1}f(x_k,\xi_k)),\\
x_{k+1}:=\Pi_X(x_k-\gamma\frac{1}{N_k}\sum^{N_k}_{j=1}f(x_{k+1/2},\xi_{k+1/2})),
\end{array}
\right.
\end{equation}
where $\gamma\in(0,\frac{1}{L})$, $L$ is the Lipschitz constant and $\Pi_X$ denotes the Euclidean projection onto $X$. The almost sure convergence for this algorithm requires the mapping to be pseudomonotone and Lipschitz continuous. However, the Lipschitz constant $L$ may be unknown or poorly estimated in practice.
To overcome this drawback, Iusem et al. \cite{I2019} presented the following variance-based extragradient with line search to solve the SVI:
\begin{equation}\label{equ:1:3}
\left\{
\begin{array}{ll}
x_{k+1/2}:=\Pi_X(x_k-\gamma_k \frac{1}{N_k}\sum^{N_k}_{j=1}f(x_k,\xi_{j,k})),\\
x_{k+1}:=\Pi_X(x_k-\gamma_k\frac{1}{N_k}\sum^{N_k}_{j=1}f(x_{k+1/2},\xi_{j,k+1/2})),
\end{array}
\right.
\end{equation}
where $\gamma_k$ is the maximum $\gamma\in\{\gamma_0\theta^{l_k}|l_k\in \mathbb{N}\cup\{0\}\}$ such that $$\gamma\|\frac{1}{N_k}\sum^{N_k}_{j=1}f(x_k+d_k,\xi_{j,k})-\frac{1}{N_k}\sum^{N_k}_{j=1}f(x_{k+1/2}(\gamma),\xi_{j,k})\|\leq\mu\|(x_k+d_k)-x_{k+1/2}(\gamma)\|$$
with $\theta\in (0,1)$ and $\mu\in (0,\frac{1}{2\sqrt{2}})$. Here,  $x_{k+1/2}(\gamma):=\Pi_X\big(x_k+d_k-\gamma(\frac{1}{N_k}\sum^{N_k}_{j=1}f(x_k,\xi_{j,k})+\beta d_k)\big)$ and $d_k$ are chosen to satisfy $\|(x_k+d_k)-x_{k+1/2}\|>0, x_k+d_k\in X,\|d_k\|\leq \sigma_k$ and $\sum\sigma_k<\infty$.
They discussed the almost sure convergence, the convergence rate and the complexity of algorithm \eqref{equ:1:3} when the mapping $F$ is H$\ddot{o}$lder/Lipschitz continuous and pseudomonotone. It is noted that Iusem et al. \cite{I2019}
pointed out that algorithm \eqref{equ:1:3} is computationally cheaper than algorithm \eqref{equ:1:2} because there is no additional information about the Lipschitz constant $L$. Morover, the iteration complexity $\mathcal{O}(1/\varepsilon)$ of algorithm \eqref{equ:1:3} is better than $\mathcal{O}(1/\varepsilon^2)$  of algorithm \eqref{equ:1:2}. Based on this pioneering work, there some variants of variance-based extragradient algorithms with line search for solving SVI have been considered in \cite{YZWL, YWL, ZDYL, LCSM}. Nevertheless, it can be observed that at each iteration the algorithms introduced in \cite{I2019, YZWL, YWL, ZDYL, LCSM} involve three oracle calls, which may not only affect the convergence speed of the algorithm but also increase the computational cost of the algorithm. Very recently,
Long and He \cite{LH} presented a fast stochastic approximated-based subgradient extragradient algorithm with line search for solving the SVI, which requires only one projection and one oracle call at each iteration, and so reducing the computation time.

\subsection{Bregman generalizations}\label{sec:1.2}
\noindent  Over the past two decades, the extragradient algorithms with the Bregman distance have attracted much attention since they have revolutionary advantages in solving optimization problems or SVIs,  see, e.g., \cite{BBC, X, GLZ, NA}).
Recently, Dang and Lan \cite{DL} presented a Bregman extragradient algorithm to solve the VI as follows:
\begin{equation*}
\left\{
\begin{array}{ll}
x_{k+1/2}:=\mathbf{P}_X(x_k,\gamma F(x_k)),\\
x_{k+1}:=\mathbf{P}_X(x_k,\gamma F(x_{k+1/2})).
\end{array}
\right.
\end{equation*}
They considered the convergence and the iteration complexity when $F$ is H$\ddot{o}$lder continuous or Lipschitz continuous.

For the SVI, Juditsky et al. \cite{JNT} proposed a stochastic Bregman extragradient algorithm (also known as stochastic mirror proximal algorithm) for solving the monotone SVI and showed a $\mathcal{O}(1/\epsilon^{2})$ oracle complexity. Chen et al. \cite{CLO} presented a accelerated stochastic Bregman extragradient algorithm for a special structure monotone SVI and obtained the $\mathcal{O}(1/\epsilon)$ iteration complexity and the $\mathcal{O}(1/\epsilon^{2})$ oracle complexity. Recently, Kannan and Shanbhag \cite{KS} introduced the following stochastic Bregman extragradient algorithm to solve the SVI:
\begin{equation}\label{equ:1:4}
\left\{
\begin{array}{ll}
x_{k+1/2}:=\mathbf{P}_X(x_k,\gamma_k f(x_k,\xi_k)),\\
x_{k+1}:=\mathbf{P}_X(x_k,\gamma_k f(x_{k+1/2},\xi_{k+1/2})),
\end{array}
\right.
\end{equation}
where $\sum^{\infty}_{k=0}\gamma^2_k<\infty$ and $\sum^{\infty}_{k=0}\gamma_k=\infty$. It is easy to see that if the distance generating function $s(x):=\frac{1}{2}\|x\|_2^2$ and the sample size $N_k=1$, then algorithm \eqref{equ:1:4} reduces to algorithm \eqref{equ:1:3}. Using the supermartingale convergence theorem, they proved the almost sure convergence of the algorithm when $F$ is pseudomonotone-plus and Lipschitz continuous. They also derived the convergence rate when $F$ is strongly pseudomonotone. However, it should be noted that there is no complexity estimates for algorithm \eqref{equ:1:4}, which is a key index to measure the performance of the SA approach. Besides, algorithm \eqref{equ:1:4} still has several points that can be improved. For instance: (i) The convergence of algorithm \eqref{equ:1:4} relies on the small stepsize policy: $\sum_k\gamma_k=\infty$ and $\sum_k\gamma^2_k<\infty$, whose salient disadvantage is that the stepsize sequence may tend to zero and thus cause dilatory convergence rate;
(ii) For every iteration $k$, algorithm \eqref{equ:1:4} only uses a single sample to approximate the value of $F$. In contrast with the multisample strategy, this may lead to high iterative complexity and slow convergence speeds;
(iii) The assumption that $F$ is pseudomonotone-plus  is somewhat strong since many functions may not satisfy this assumption.

Based on the above, a natural question is: Can we establish an improved stochastic Bregman extragradient algorithm with multisamples, an unknown Lipschitz constant and perfect complexity results to solve the SVI when $F$ is not necessarily pseudomonotone?  Our aim in this paper is to answer the above question in the affirmative.

\subsection{Contributions}\label{sec:1.3}
\noindent Motivated by the works of \cite{I2019, DL, KS}, the overarching goal of this paper is to propose a variance-based stochastic Bregman extragradient algorithm with line search for solving the SVI.
Our contributions are summarized as follows:

\begin{itemize}
\item  We employ a novel stochastic Bregman line search rule to fix the stepsize which makes the stepsize away from zero and does not require the information of the Lipschitz constant. Moreover, our algorithm only needs to compute a single oracle call at each iteration when $X$ is unbounded. These may improve the efficiency and applicability of the proposed algorithm.

\item We develop a more concise and effective way to analyze the almost sure convergence of the algorithm without using the supermartingale convergence theorem, which was used in \cite{I2019, YZWL, DL, KS}. Moreover, we do not require that $F$ is pseudomonotone-plus used in \cite{KS}.

\item Under the Minty variational inequality condition, we establish the convergence rate in terms of the gap function and the natural residual function, respectively. Moreover, the iteration complexity and the oracle complexity of our algorithm are discussed.

\item  We give some numerical experiments of our algorithm and compare it with the algorithms considered in \cite{I2019,YZWL,KS}. The numerical results show that our algorithm is competitive with other algorithms, and then we apply it to solve the stochasic Nash game problem.
\end{itemize}

The rest of the paper is organized as follows: In Section \ref{sec:Preliminaries}, we review some related definitions and results. The new algorithm and its convergence analysis are presented in Section \ref{sec:Main}. In Section \ref{sec:Convergence rate and complexity}, we discuss the convergence rate and the complexity of the algorithm under two different cases. Several numerical results are given in Section \ref{sec:Numerical} to show the competitiveness of the proposed algorithm and the conclusions are discussed in Section \ref{sec:conclusion}.

\section{Preliminaries}\label{sec:Preliminaries}
\setcounter{equation}{0}

\noindent In this section, we recall some notations and results that will be used later on. For $x,y\in \mathbb{R}^n$, we denote by $\left\|x\right\|$ the Euclidean norm, and $\|x\|_*=\sup\{z^Tx|\|z\|\leq 1\}$ the dual norm. For a random variable $\xi$ and a $\sigma$-algebra $\mathcal{F}$, we denote by $\mathbb{E}[\xi]$ the expectation of $\xi$, $\mathbb{E}[\xi|\mathcal{F}]$ the conditional expectation of $\xi$ with respect to $\mathcal{F}$, $V(\xi)$ the variance of $\xi$, and $V(\xi|\mathcal{F})$ the conditional variance of $\xi$ with respect to $\mathcal{F}$.
For $p\geq1$, $|\xi|_p$ is the $\mathcal{L}_p$-norm of $\xi$ and $|\xi|\mathcal{F}|_p:=\sqrt[p]{\mathbb{E}[|\xi|_p|\mathcal{F}]}$ is the $\mathcal{L}_p$-norm of $\xi$ conditional to $\mathcal{F}$.
We write $\xi\in \mathcal{F}$ if $\xi$ is $\mathcal{F}$-measurable, $\xi\perp\perp\mathcal{F}$ if $\xi$ is
independent of $\mathcal{F}$, i.i.d. for independently identical distributed, and a.s. for `` almost surely".
Let $\mathbb{N}_0:=\mathbb{N}\cup\{0\}$, where $\mathbb{N}$ is the set of all positive integers. For a given number $p\geq 2$, we define the oracle error $\epsilon:\mathbb{R}^n\times\Xi\rightarrow\mathbb{R}$ for \eqref{equ:1:1} by
$$\epsilon(x,\xi):=f(x,\xi)-F(x),\; \forall \; \xi\in \Xi,x\in \mathbb{R}^n.$$

Given i.i.d samples $\xi_N:=\{\xi\}^N_{j=1}$ drawn from $\Xi$, the empirical mean operator is defined by
$$\widehat F(x;\xi_N):=\frac{1}{N}\sum_{j=1}^{N}f(x,\xi_{j}),~~\forall x\in\mathbb{R}^n.$$
A function $s:X\subseteq\mathbb{R}^n\to \mathbb{R}$ is called a distance generating function, if $s$ is continuously differentiable and $\alpha$-strongly convex on $X$ ($\alpha>1$), i.e.,
\begin{equation*}
\big(\nabla{s(x)}-\nabla{s(y)}\big)^T(x-y)\geq \alpha\|x-y\|^2,\;\forall x, y\in X.
\end{equation*}
The Bregman distance $V(x,z)$ associated with $s$ is defined by
\begin{equation}\label{equ:2:1}
V(x,z):=s(z)-\Big[s(x)+\nabla{s(x)}^T(z-x)\Big],
\end{equation}
which was initially studied by Bregman \cite{B} and later by many others (see, e.g., \cite{AT,BBC,T}). By \eqref{equ:2:1} and the strong convexity of $s$, we can easily get
\begin{equation}\label{equ:2:2}
V(x,z)\geq\frac{\alpha}{2}\|x-z\|^2,\;\forall x,z\in X.
\end{equation}
Note that $V(x,z)$ reduces to the Euclidean distance when $s(x)=\frac{1}{2}\|x\|^2$. Furthermore, if the gradient of the distance generating function is $\mathcal{Q}$-Lipschitz continuous, i.e., there exists $\mathcal{Q}>0$ such that
$$\|\nabla{s(x)}-\nabla{s(z)}\|_*\leq \mathcal{Q}\|x-z\|,\;\forall x,z\in X,$$
then (see, Lemma 1.2.3 of \cite{N})
\begin{equation}\label{equ:2:3}
V(x,z)\leq \frac{\mathcal{Q}}{2}\|x-z\|^2,\;\forall x,z\in X.
\end{equation}
We say that $V$  is growing quadratically when \eqref{equ:2:3} holds.

For $x\in X$ and $r\in \mathbb{R}^n$, the resulting prox-mapping $\mathbf{P}_{X}: X\times \mathbb{R}^n\to X$ is defined by \begin{equation}\label{equ:2:4}
\mathbf{P}_{X}(x,r):=\mathop{\mathrm{argmin}}\limits_{z\in X}\Big(r^Tz+V(x,z)\Big).
\end{equation}
Moreover, for a given operator $H:\mathbb{R}^n\to \mathbb{R}^n$ and a constant $a>0$, we define the natural residual function $R_a$ by
\begin{equation}\label{equ:2:5}
R_{a}(x,H):=\|x-\mathbf{P}_X(x,aH(x))\|,\;\forall x\in \mathbb{R}^n.
\end{equation}

We need the following importance assumptions.

\begin{assumption}\label{ass:2.1}
\begin{itemize}
\item[(i)] There exists a measurable function $L:\Xi\times \Xi\to \mathbb{R}_+$ such that $L(\xi,\eta)\geq 1$ for almost every $\xi, \eta \in \Xi$ and $$\left\|f(x,\xi)-f(y,\eta)\right\|\leq L(\xi,\eta)\left\|x-y\right\|,\; \forall x,y\in \mathbb{R}^n.$$

\item[(ii)] There exist $x^*\in \mathbb{R}^n$ and $p\geq 2$ such that $\mathbb{E}[\left\|f(x^*,\xi)\right\|^p]<\infty $ and $\mathbb{E}_{(\xi,\eta)}[L(\xi,\eta)^p]<\infty .$ In what follows, we denote by $\mathcal{L}:=\mathbb{E}_{(\xi,\eta)}[L(\xi,\eta)]$ and $\mathcal{L}_q:=|L(\xi,\eta)|_q+\mathcal{L}$ for any $q>0.$
\end{itemize}
\end{assumption}
%\begin{assumption}
%\begin{itemize}
%\item[(i)] There exists a measurable function $L:\Xi \to \mathbb{R}_+$ such that $ L(\xi)\geq 1$ for almost every $\xi \in \Xi$ and $$\left\|f(x,\xi)-f(y,\xi)\right\|\leq L(\xi)\left\|x-y\right\|,\; \forall x,y\in \mathbb{R}^n.$$
%
%\item[(ii)] There exist $x^*\in \mathbb{R}^n$ and $p\geq 2$ such that $\mathbb{E}[\left\|f(x^*,\xi)\right\|^p]<\infty $ and $\mathbb{E}[ L(\xi)^p]<\infty .$ In what follows, we denote by $\mathcal{L}:=\mathbb{E}[L(\xi)]$ and $\mathcal{L}_q:=|L(\xi)|_q+\mathcal{L}$ for any $q>0.$
%\end{itemize}
%\end{assumption}

The following lemma is similar to Lemma 1.2 in \cite{I2019}, which is a simple consequence of Jensen and Minkowski's inequalities, and so we omit its proof.

\begin{lemma}\label{lemma:2.1}
Let Assumption \ref{ass:2.1} hold. Then for any $q\in[p,2p]$, the mapping $F$ is $\mathcal{L}$-Lipschitz continuous on $\mathbb{R}^n$ and $|\|\epsilon(\cdot,\xi)\||_q$ is $\mathcal{L}_q$-Lipschitz continuous on $\mathbb{R}^n.$
\end{lemma}

We also need the following lemmas.

\begin{lemma} \label{lemma:2.2}  \cite{AL}
Given a convex function $f$, for any $x_1, \cdots, x_n\in \mathrm{dom}f$ and constants $\alpha_1, \cdots, \alpha_n$ satisfying $\sum^{n}_{i=1}\alpha_i=1, \alpha_i\geq0, i=1,\cdots, n,$ we have $f(\sum^{n}_{i=1}\alpha_ix_i)\leq\sum^{n}_{i=1}\alpha_if(x_i).$ In other word, $f(\mathbb{E}[x])\leq\mathbb{E}[f(x)].$
\end{lemma}

\begin{lemma}\label{lemma:2.3} \cite{LGH}
Let $X\in \mathbb{R}^n$ be a convex set and $p:X\to \mathbb{R}$ be a differentiable convex function. If $u^*$ is an optimal solution of $\min\{p(u)+V(\tilde{x},u):u\in X\}$, then
$$p(u^*)+V(u^*,u)+V(\tilde {x},u^*)\leq p(u)+V(\tilde {x},u),~~ \forall u\in X.$$
\end{lemma}

\section{Algorithm and convergence analysis}\label{sec:Main}
\setcounter{equation}{0}

\noindent In this section, we propose an improved stochastic Bregman extragradient algorithm with line search to solve the SVI and give its convergence analysis.

\subsection{Algorithm }\label{sec:Main1}
\noindent  Inspired by the works of \cite{KS,DL}, we proposed the following stochastic Bregman extragradient algorithm with line search for solving the SVI.

\begin{table*}[h]
\fbox{%
\parbox{\textwidth}{%
%\noindent\rule[0.25\baselineskip]{\textwidth}{1pt}
\noindent{\bf Algorithm 1 (Stochastic Bregman extragradient algorithm with line search)} \\
%\noindent\rule[0.25\baselineskip]{\textwidth}{1pt}
\noindent{\bf Step 0:} Choose an initial point $x_0\in X, \gamma_0,\theta \in (0,1) $, and the sample rate $\{N_{k}\}\subset \mathbb{N}$. Set $k:=0$.

\noindent{\bf Step 1:}  Generate a sample $\xi_k:=\{\xi_{j,k}\}_{j=1}^{N_{k}}$ from $\Xi$ and compute $\widehat{F}(x_k;\xi_{k})=\frac{1}{N_{k}}\sum_{j=1}^{N_{k}}f(x_k,\xi_{j,k})$. If
$x_k=\mathbf{P}_{X}(x_k,\theta^{-1}\gamma_0\widehat{F}(x_k;\xi_{k}))$, regenerate a sample. Otherwise, go to \noindent{\bf Step 2}.

\noindent{\bf Step 2:}  Generate a sample $\xi_{j,k+1/2}:=\{\xi_{j,k+1/2}\}^{N_k}_{j=1}$ from $\Xi$ and choose $\gamma_{k}$ as the maximum $\gamma \in \{\gamma_0 \theta^{l_k}|l_k\in \mathbb{N}_{0}\}$ such that \begin{equation}\label{equ:3:1}
\gamma^2\|\widehat{F}(x_k;\xi_{k})-\widehat{F}(x_{k+1/2}(\gamma);\xi_{k+1/2})\|^2\leq \alpha V(x_k,x_{k+1/2}(\gamma )),
\end{equation}
where $$x_{k+1/2}(\gamma)=\mathbf{P}_{X}(x_k,\gamma\widehat{F}(x_k;\xi_{k})),\;\widehat{F}(x_{k+1/2}(\gamma );\xi_{k+1/2})=\frac{1}{N_{k}}\sum_{j=1}^{N_{k}}f(x_{k+1/2}(\gamma),\xi_{j,k+1/2}).$$

\noindent{\bf Step 3:} Compute
\begin{equation}\label{equ:3:2}
x_{k+1/2}=\mathbf{P}_{X}(x_k,\gamma_k\widehat{F}(x_k;\xi_{k})),
\end{equation}
\begin{equation}\label{equ:3:3}
x_{k+1}=\mathbf{P}_{X}(x_k,\gamma_k\widehat{F}(x_{k+1/2};\xi_{k+1/2})).
\end{equation}
Set $k:=k+1$ and return to \noindent{\bf Step 1}.
}%
}
\end{table*}

\begin{remark}
\begin{itemize}
\item [(i)]  Compared with the algorithm in\cite{KS}, we give a stochastic Bregman version of Armijo's line search rule to select the stepsize at each iteration, which not only is robust with respect to unknown Lipschitz constant but also allows us to build an adaptive stepsize bounded away from zero. Moreover, Algorithm 1 only involves one oracle error when discussing the convergence rate of the algorithm in the case that $X$ is unbounded.

\item [(ii)] Different from the single sample strategy used in \cite{KS}, we adopt the multisample in Algorithm 1 to approximate the value of $F$ at each iteration, which may bring faster iteration complexity and more superiority in dealing with large-scale problems.

\item [(iii)] In Algorithm 1, since we adopt a sample-based line search rule, the termination criteria is generally unclear. Hence, we choose the strategy of regenerating a sample in Step 1, which can ensure the line search is well-defined and the algorithm terminates in a finite number of iterations.
\end{itemize}
\end{remark}

In the following, let the solution set $X^*$ of the SVI be nonempty. We start with some necessary assumptions that will be used later on.

\begin{assumption}\label{ass:3.1}
 \begin{itemize}
\item [(i)] For any $x\in X,$ there exists $B_*$ such that $\|F(x)\|_*\leq \frac{B_*}{2}.$

\item [(ii)] For any $x\in X,$ there exists a constant $v>0$ such that $|\|f(x,\xi)-F(x)\|_*|_2\leq v.$

\end{itemize}
\end{assumption}

\begin{assumption}\label{ass:3.2}
 For any $y\in X^*,$ the following Minty variational inequality holds:
$$F(x)^T(x-y)\geq0, ~~\forall x\in X.$$
\end{assumption}

\begin{remark}
We observe that Assumption \ref{ass:3.2} holds if $F$ is pseudomonotone. However, the converse
does not hold; see \cite{DL} for a counter-example. Therefore, Assumption \ref{ass:3.2} is weaker than the pseudomonotonity considered in \cite{I2019}.
\end{remark}

\begin{assumption}\label{ass:3.3}
In Algorithm $1$, the sequences $\{\xi_{j,k}|k\in \mathbb{N}_0, j\in [N_k]\}$ and $\{\xi_{j,k+1/2}|k\in \mathbb{N}_0, j\in [N_k]\}$ are i.i.d. samples drawn from $\Xi$ and independent of each other. Moreover, $\sum_{k=0}^{\infty}\frac{1}{N_k}<\infty $.
\end{assumption}

Let $\mathcal{F}_k:=\sigma\{x_0,\xi_0,\xi_{1/2},\xi_{1},\cdots,\xi_{k-1/2}\}$ and $\mathcal{F}_{k+1/2}:=\mathcal{F}_k\cup\{\xi_k\}$ be the $\sigma$-algebra which relate to the generation of $x_k, x_{k+1/2}$ and $\gamma_k$.
The oracle errors are defined as follows: $$\epsilon^k_1:=\widehat {F}(x_k;\xi_k)-F(x_k),~~ \epsilon^k_2:=\widehat {F}(x_{k+1/2};\xi_{k+1/2})-F(x_{k+1/2}),~~\epsilon^k_3:=\widehat {F}(x_{k+1/2};\xi_{k})-F(x_{k+1/2}).$$
It is easy to see that $x_k\in\mathcal{F}_k, \;x_{k+1/2}\in\mathcal{F}_{k+1/2}$ and
\begin{equation}\label{equ:3:4}
\mathbb{E}[\epsilon^k_1|\mathcal{F}_k]=0,\; \mathbb{E}[\epsilon^k_2|\mathcal{F}_k]=\mathbb{E}[\mathbb{E}[\epsilon^k_2|\mathcal{F}_{k+1/2}]|\mathcal{F}_k]=0,\; \mathbb{E}[\epsilon^k_3|\mathcal{F}_k]\neq0.
\end{equation}

The following lemma provides a simple characterization of a solution to the SVI.

\begin{lemma} \label{lemma:3.1}
A point $x\in X$ is a solution of the SVI if and only if
\begin{equation*}
x=\mathbf{P}_X(x,\gamma \widehat F(x;\xi_N)),
\end{equation*}
where $F(x)+\epsilon=\widehat F(x;\xi_N):=\frac{1}{N}\sum^{N}_{j=1}f(x,\xi_j)$ and $\gamma>0.$
\end{lemma}

\noindent{\bf Proof.} If $x=\mathbf{P}_X(x,\gamma \widehat F(x;\xi_N))$, then by the optimality condition of \eqref{equ:2:4},
\begin{equation*}
\Big(\gamma \widehat F(x;\xi_N)+\nabla {s(\mathbf{P}_X(x,\gamma \widehat F(x;\xi_N)))-\nabla {s(x)}\Big)^T\Big(z-\mathbf{P}_X(x,\gamma \widehat F(x;\xi_N))}\Big)\geq 0,~~\forall z\in X.
\end{equation*}
From $\widehat{F}(x;\xi_N)=F(x)+\epsilon$, we have $\gamma (F(x)+\epsilon)^T(z-x)\geq0$ for any $z\in X$. If $x\in \mathcal{F}$ and $\mathbb{E}[\epsilon|\mathcal{F}]=0$, by taking $\mathbb{E}[\cdot|\mathcal{F}]$ in the above inequality, we obtain$$\gamma\mathbb{E}[F(x)|\mathcal{F}]^T(z-x)+\gamma\mathbb{E}[\epsilon|\mathcal{F}]^T(z-x)\geq0,$$
which equals to $$\gamma F(x)^T(z-x)\geq0.$$
Since $\gamma>0$, $x$ is a solution of the SVI. The ``only if " part of the statement easily follows from the optimality condition of \eqref{equ:2:4}. This completes the proof.$\hfill\square$

By \eqref{equ:2:5} and Lemma \ref{lemma:3.1}, we have the following corollary.

\begin{cor}\label{cor:3.1}
A point $x\in X$ is a solution of the SVI if and only if $\|R_a(x,H)\|=0$ for a constant $a>0.$
\end{cor}

The following proposition shows that the line search rule is well-defined under Assumption \ref{ass:2.1}.

\begin{proposition}
Let Assumption \ref{ass:2.1} hold. Then the line search rule \eqref{equ:3:1} terminates in a finite number of iterations.
\end{proposition}
\noindent{\bf Proof.}
We consider the following two cases. (i) If $x_k\in X^*$, then we obtain by Lemma \ref{lemma:3.1} with $H:=\frac{1}{N_k}\sum^{N_k}_{j=1}f(x_k,\xi_{j,k})$ that
\begin{equation}\label{equ:3:5}
\|x_k-x_{k+1/2}(\gamma)\|^2=\|x_k-\mathbf{P}_X(x_k,\gamma\widehat {F}(x_k;\xi_k))\|^2=\|R_\gamma(x_k,H)\|^2.
\end{equation}
By Corollary \ref{cor:3.1}, $\|R_\gamma(x_k,H)\|=0$. It follows from \eqref{equ:2:3}, \eqref{equ:3:1} and \eqref{equ:3:5} that
$$\gamma^2\|\widehat{F}(x_k;\xi_{k})-\widehat{F}(x_{k+1/2}(\gamma);\xi_{k+1/2})\|^2\leq \alpha V(x_k,x_{k+1/2}(\gamma))\leq\frac{\alpha\mathcal{Q}}{2}\|x_k-x_{k+1/2}(\gamma)\|^2=0.$$
Since $\lim\limits_{k\to\infty}\gamma_0\theta^{l_k}=0$, taking the limit $k\to \infty$ in above inequality, we get
$$0\leq \alpha V(x_k,x_{k+1/2}(\gamma))\leq0.$$
Hence, we can get the conclusion immediately.

In the case of $x_k\notin X^*$ but $x_k\in X$,  suppose by contradiction that the line search step does not terminate in a finite number of iterations, which means that, for every $k$,
\begin{equation*}
\gamma^2\|\widehat{F}(x_k;\xi_{k})-\widehat{F}(x_{k+1/2}(\gamma);\xi_{k+1/2})\|^2> \alpha V(x_k,x_{k+1/2}(\gamma)),
\end{equation*}
which together with \eqref{equ:2:2} implies
\begin{equation}\label{equ:3:6}
\gamma^2\|\widehat{F}(x_k;\xi_{k})-\widehat{F}(x_{k+1/2}(\gamma);\xi_{k+1/2})\|^2>\frac{\alpha^2}{2}\|x_k-x_{k+1/2}(\gamma)\|^2.
\end{equation}
Taking the limit $k\to \infty$ in \eqref{equ:3:6}, we obtain $0\geq\frac{\alpha^2}{2}\|x_k-x_{k+1/2}(\gamma)\|^2\|>0,$ which gives a contradiction.

(ii) If $x^k\notin X$, then we suppose that the line search step does not terminate in a finite number of iterations, by taking the limit $k\to \infty$ in \eqref{equ:3:6}, we have for every $k$ that $$0\geq\frac{\alpha^2}{2}\|x_k-\mathbf{P}_X(x_k,\gamma\hat {F}(x_k,\xi_{k}))\|^2>0,$$ which is a contradiction. Therefore, the line search Sep 2 terminates in a finite number of iterations. This completes the proof.$\hfill\square$

\subsection{Convergence analysis}\label{sec:Main2}
\noindent In this subsection, we aim to discuss the a.s. convergence of Algorithm 1. Since we use a stochastic line search rule to determine the stepsize in Algorithm 1, we need to control the lower bound of $|\gamma_k|\mathcal{F}_k|_2$.

\begin{lemma}\label{lemma:3.2}
Let $\bar {L}(\xi_k,\xi_{k+1/2}):=\frac{1}{N_k}\sum_{j=1}^{N_k}L(\xi_{j,k},\xi_{j,k+1/2})$ and Assumptions \ref{ass:2.1} and \ref{ass:3.3} hold. Then we have a.s. $\gamma_k\geq \min\{\frac{\alpha\theta}{\bar{L}(\xi_k,\xi_{k+1/2})\sqrt{2}}, \gamma_0\}$ and $|\gamma_k|\mathcal{F}_k|_2\cdot |L(\xi,\xi_{+1/2})|_2\geq \min \{\frac{\alpha\theta}{\sqrt{2}},\gamma_0 \}.$
\end{lemma}

\noindent{\bf Proof.} If $\gamma_0$ satisfies \eqref{equ:3:1}, then $\gamma_k=\gamma_0.$ Otherwise, one has
$$\theta^{-2}\gamma^2_k\|\widehat{F}(x_k;\xi_{k})-\widehat{F}(x_{k+1/2}(\theta^{-1}\gamma_k);\xi_{k+1/2})\|^2>\alpha V(x_k,x_{k+1/2}(\theta^{-1}\gamma_k)).$$
It follows from \eqref{equ:2:2} that
\begin{equation}\label{equ:3:7}
\theta^{-2}\gamma^2_k\|\widehat{F}(x_k;\xi_{k})-\widehat{F}(x_{k+1/2}(\theta^{-1}\gamma_k);\xi_{k+1/2})\|^2>\frac{\alpha^2}{2}\|x_k-x_{k+1/2}(\theta^{-1}\gamma_k)\|^2.
\end{equation}
By the definition of $\widehat{F}(x;\xi_N)$ and Assumption \ref{ass:2.1},
\begin{equation}\label{equ:3:8}
\|\widehat{F}(x_k;\xi_{k})-\widehat{F}(x_{k+1/2}(\theta^{-1}\gamma_k);\xi_{k+1/2})\|^2\leq\bar L(\xi_k,\xi_{k+1/2})^2\|x_k-x_{k+1/2}(\theta^{-1}\gamma_k)\|^2.
\end{equation}
Combining \eqref{equ:3:7}, \eqref{equ:3:8} and the fact that $x_k\neq x_{k+1/2}(\theta^{-1}\gamma_k)$, we get
$$\gamma^2_k\geq\frac{\alpha^2\theta^2}{2\bar{L}(\xi_k,\xi_{k+1/2})^2}.$$
Since $\gamma_k>0,$ it follows that
$$\gamma_k\geq\frac{\alpha\theta}{\bar{L}(\xi_k,\xi_{k+1/2})\sqrt{2}}.$$

On the other hand, since a.s. $L(\xi_k,\xi_{k+1/2})\geq1,$ we have a.s. $\gamma_k\bar{L}(\xi_k,\xi_{k+1/2})\geq\min\{\frac{\alpha\theta}{\sqrt{2}},\gamma_0\}.$
By Assumption \ref{ass:3.3},
\begin{align*}
\min\{\frac{\alpha\theta}{\sqrt{2}},\gamma_0\}\leq&\mathbb{E}[\gamma_k\bar{L}(\xi_k,\xi_{k+1/2})|\mathcal{F}_k]\\
\leq&|\gamma_k|\mathcal{F}_k|_2\cdot|\bar{L}(\xi_k,\xi_{k+1/2})|\mathcal{F}_k|_2\\
\leq&|\gamma_k|\mathcal{F}_k|_2\sqrt{\frac{1}{N_k}\sum^{N_k}_{j=1}\mathbb{E}[L(\xi_{j,k},\xi_{j,k+1/2})^2|\mathcal{F}_k]}\\
=&|\gamma_k|\mathcal{F}_k|_2\cdot|L(\xi,\xi_{+1/2})|_2.
\end{align*}
This completes the proof.$\hfill\square$

\begin{remark} By Lemma \ref{lemma:3.2}, the stepsize sequence $\{\gamma_k\}$ generated by Algorithm $1$ satisfies $|\gamma_k|\mathcal{F}_k|_2\cdot |L(\xi,\xi_{+1/2})|_2\geq \min \{\frac{\alpha \theta}{\sqrt{2}},\gamma_0 \}$ and $\gamma_k<1$ for all $k$; that is, the conditional expectation to $\mathcal{F}_k$ of the sequence $\{\gamma_k\}$ is bounded.
\end{remark}

In \cite{I2019, KS}, the authors obtained the a.s. convergence of the algorithm by using the supermartingale convergence theorem. This requires that a strict recursive relation must be satisfied.  However, many random variable sequences may not satisfy this relation. Therefore, we will introduce a more concise and efficient way to establish the a.s. convergence as follows.

\begin{lemma}\label{lemma:3.3}
Let $x_k\in X$ and $\gamma_k\in(0,1), k=0, 1, \cdots$. Assume that the following conditions hold:
\begin{itemize}
\item [(i)] $\lim\limits_{k\to\infty}V(x_k,\mathbf{P}_X(x_k,\gamma_k\widehat F(x_k;\xi_k)))=0;$

\item [(ii)] there exists $K\in\mathbb{N}$ and $\gamma^*>0$ such that $\gamma_k\geq \gamma^*$ for any $k\geq K.$
\end{itemize}
Then we have a.s. $\lim\limits_{k\to \infty}\|R_{\gamma_k}(x_k,H)\|=0.$
If in addition, the sequence $\{x_k\}$ is bounded, there a.s. exists an accumulation point $\hat x$ of $\{x_k\}$ such that $\hat x\in X^*$.
\end{lemma}

\noindent{\bf Proof.} Set $x_{k+1/2}=\mathbf{P}_{X}(x_k,\gamma_k\widehat{F}(x_k;\xi_{k})).$  It follows from  (i) and \eqref{equ:2:2} that $\lim\limits_{k\to\infty}\|x_k-x_{k+1/2}\|=0.$
By (ii) and the definition of $R_{a}(x,H)$ in \eqref{equ:2:5} with
$H:=\frac{1}{N_k}\sum^{N_k}_{j=1}f(x_k,\xi_{j,k})$, we have a.s. $$\lim_{k\to\infty}\|R_{\gamma_k}(x_k,H)\|=0.$$
Moreover, if the sequence $\{x_k\}$ is bounded, there exists a subsequence $\{\hat {x}_i\}$ of $\{x_k\}$ obtained by setting $\hat {x}_i=x_{n_i}$ for $n_1\leq n_2\leq \cdots,$ such that $\lim\limits_{i\to \infty} \|\hat {x}_i-\hat x\|=0.$ Let $\{\hat{x}_{i+1/2}\}$ be the corresponding subsequence in $\{x_{k+1/2}\}$, i.e., $x_{i+1/2}=\mathbf{P}_X(x_{n_i},\gamma_{n_i}\widehat{F}(x_{n_i};\xi_{n_i}))$ and $\hat \gamma_i=\gamma_{n_i}$, it follows that $\lim\limits_{i\to \infty}\|\hat {x}_i-\hat {x}_{i+1/2}\|=0.$
Furthermore, by the optimality condition of \eqref{equ:2:4},
$$\Big(\hat F(\hat x_i;\xi_i)+\frac{1}{\hat {\gamma_i}}\big[\nabla{s}(\hat{x}_{i+1/2})-\nabla{s}(\hat{x}_i)\big]\Big)^T\Big(z-\hat {x}_{i+1/2}\Big)\geq0,\forall z\in X, i\geq 1.$$
Taking the limit $i\to +\infty$ in the above inequality and using the continuity of $H(\cdot)$ and $\nabla s(\cdot)$, we have a.s. from condition (ii) that $\hat F(\hat x_i;\xi_i)^T(z-\hat {x}_{i+1/2})\geq0$ for any $z\in X.$ By taking $\mathbb{E}[\cdot|\mathcal{F}_k]$ in it, we get a.s. $F(\hat x)^T(z-\hat x)\geq0$ for any $z\in X.$ This completes the proof.
$\hfill\square$

\begin{remark}
It is worth mentioning that Lemma \ref{lemma:3.3} provides a unified framework to analyse the a.s. convergence of the algorithm, which only relies on the natural residual function and the lower bound of the stepsize sequence. Actually, Lemma \ref{lemma:3.3} can be reviewed as an extension of the stochastic case of Lemma 3 in \cite{DL} and so it covers a wider range of applications.
\end{remark}

We now show the a.s. convergence of Algorithm 1.

\begin{theorem}\label{the:3.1}
Let Assumption \ref{ass:3.2} hold. Then for any $x^*\in X^*,$ the sequence $\{x_k\}$ and $\{x_{k+1/2}\}$ generated by Algorithm 1 are bounded a.s. and converge a.s. to a point in $X^*$.
\end{theorem}

\noindent{\bf Proof.} By \eqref{equ:3:2} and Lemma \ref{lemma:2.3},
$$\gamma_k\widehat {F}(x_k;\xi_k)^T(x_{k+1/2}-x)+V(x_{k+1/2},x)+V(x_k,x_{k+1/2})\leq V(x_k,x),~~\forall x\in X.$$
Setting $x=x_{k+1}$ in the above inequality, we get
\begin{equation}\label{equ:3:9}
\gamma_k\widehat  {F}(x_k;\xi_k)^T(x_{k+1/2}-x_{k+1})+V(x_{k+1/2},x_{k+1})+V(x_k,x_{k+1/2})\leq V(x_k,x_{k+1}).
\end{equation}
Similarly, from \eqref{equ:3:3} and Lemma \ref{lemma:2.3}, we derive
\begin{equation}\label{equ:3:10}
\gamma_k\widehat  {F}(x_{k+1/2};\xi_{k+1/2})^T(x_{k+1}-x)+V(x_{k+1},x)+V(x_k,x_{k+1})\leq V(x_k,x),~~\forall x\in X.
\end{equation}
Noting that $\widehat {F}(x_{k+1/2};\xi_{k+1/2})^T(x_{k+1}-x)=\widehat  {F}(x_{k+1/2};\xi_{k+1/2})^T(x_{k+1/2}-x)-\widehat {F}(x_{k+1/2};\xi_{k+1/2})$
$^T(x_{k+1/2}-x_{k+1}),$ letting $x=x^*$ in \eqref{equ:3:10} and substituting \eqref{equ:3:9} into \eqref{equ:3:10}, we obtain
\begin{align*}
 V(x_k,x^*)\geq&\gamma_k\widehat {F}(x_{k+1/2};\xi_{k+1/2})^T(x_{k+1/2}-x^*)-\gamma_k\widehat  {F}(x_{k+1/2};\xi_{k+1/2})^T(x_{k+1/2}-x_{k+1})\\&+V(x_{k+1},x^*)
+\gamma_k\widehat  {F}(x_k;\xi_k)^T(x_{k+1/2}-x_{k+1})+V(x_{k+1/2},x_{k+1})+V(x_k,x_{k+1/2}).
\end{align*}
This implies that
\begin{align*}
 V(x_k,x^*)\geq&\gamma_k(\widehat {F}(x_k;\xi_k)-\widehat {F}(x_{k+1/2};\xi_{k+1/2}))^T(x_{k+1/2}-x_{k+1})+
V(x_{k+1},x^*)\\ &+\gamma_k(\epsilon^k_2+F(x_{k+1/2}))^T(x_{k+1/2}-x^*)+V(x_{k+1/2},x_{k+1})+V(x_k,x_{k+1/2}).
\end{align*}
From Assumption \ref{ass:3.2}, we have $F(x_{k+1/2})^T(x_{k+1/2}-x^*)\geq 0$. It follows that
\begin{align}\label{equ:3:11}
V(x_k,x^*)\geq&\gamma_k(x_{k+1/2}-x^*)^T\epsilon^k_2+\gamma_k(\widehat{F}(x_k;\xi_k)-\widehat {F}(x_{k+1/2};\xi_{k+1/2}))^T(x_{k+1/2}-x_{k+1})\nonumber\\
&+V(x_{k+1},x^*)+V(x_{k+1/2},x_{k+1})+V(x_k,x_{k+1/2}).
\end{align}
On the other hand, we have from \eqref{equ:2:2} that
\begin{align*}
\gamma_k(\widehat  {F}(x_k;\xi_k)-&\widehat {F}(x_{k+1/2};\xi_{k+1/2}))^T(x_{k+1/2}-x_{k+1})+V(x_{k+1/2},x_{k+1})+V(x_{k},x_{k+1/2})\\
\geq& -\gamma_k\|\widehat {F}(x_k;\xi_k)-\widehat  {F}(x_{k+1/2};\xi_{k+1/2})\|\cdot \|x_{k+1/2}-x_{k+1}\|\\
&+V(x_{k+1/2},x_{k+1})+V(x_{k},x_{k+1/2})\\
\geq& -\gamma_k\|\widehat  {F}(x_k;\xi_k)-\widehat  {F}(x_{k+1/2};\xi_{k+1/2})\|\cdot\sqrt{\frac{2}{\alpha}V(x_{k+1/2},x_{k+1})}\\
&+V(x_{k+1/2},x_{k+1})+V(x_{k},x_{k+1/2})\\
\geq &-\frac{\gamma_k^2}{2\alpha}\|\widehat  {F}(x_k;\xi_k)-\widehat {F}(x_{k+1/2};\xi_{k+1/2})\|^2+V(x_k,x_{k+1/2}).
\end{align*}
It follows from \eqref{equ:3:11} that
$$\gamma_k(x_{k+1/2}-x^*)^T\epsilon^k_2-\frac{\gamma_k^2}{2\alpha}\|\widehat {F}(x_k;\xi_k)-\widehat {F}(x_{k+1/2};\xi_{k+1/2})\|^2+V(x_k,x_{k+1/2})\leq V(x_k,x^*)-V(x_{k+1},x^*).$$
This together with \eqref{equ:3:1} yields
$$\gamma_k(x_{k+1/2}-x^*)^T\epsilon^k_2+\frac{1}{2}V(x_k,x_{k+1/2})\leq V(x_k,x^*)-V(x_{k+1},x^*).$$
Taking $\mathbb{E}[\cdot|\mathcal{F}_k]$ in the above inequality, we have from \eqref{equ:3:4} that
\begin{equation}\label{equ:3:12}
\frac{1}{2}\mathbb{E}[V(x_k,x_{k+1/2})|\mathcal{F}_k]\leq \mathbb{E}[V(x_k,x^*)|\mathcal{F}_k]-\mathbb{E}[V(x_{k+1},x^*)|\mathcal{F}_k].
\end{equation}
It is easy to see that the sequence $\mathbb{E}[V(x_k,x^*)|\mathcal{F}_k]$ is nonincreasing. Since $x_k\in \mathcal{F}_k,$ it follows that the sequence $\{V(x_k,x^*)\}$ is nonincreasing. Therefore, $\{V(x_k,x^*)\}$ is a convergent sequence and the sequence $\{x_k\}$ is bounded in the sense of almost surely.

Moreover, summing up the inequality in \eqref{equ:3:12} from $k=0$ to $\infty$, we obtain
$$\frac{1}{2}\sum^{\infty}_{k=0}\mathbb{E}[V(x_k,x_{k+1/2})|\mathcal{F}_k]\leq\mathbb{E}[V(x_0,x^*)|\mathcal{F}_k].$$
Taking the expectation in the above inequality, we get
$$\frac{1}{2}\sum^{\infty}_{k=0}\mathbb{E}[V(x_k,x_{k+1/2})]\leq \mathbb{E}[V(x_0,x^*)]<\infty,$$
which implies that $$\lim\limits_{k\to +\infty}\mathbb{E}[V(x_k,x_{k+1/2})]=0.$$
Since the sequences $\{V(x_k,x_{k+1/2})\}$ are discrete variables, and for any $k\geq0,\;V(x_k,x_{k+1/2})\geq0$, thus
\begin{equation}\label{equ:3:13}
\lim\limits_{k\to +\infty}V(x_k,x_{k+1/2})=0.
\end{equation}
By invoking Lemmas \ref{lemma:3.2}-\ref{lemma:3.3} and the definition of $\{x_{k+1/2}\}$, there exists an accumulation point $\hat x$ of $\{x_k\}$ such that $\hat x\in X^*.$ We can replace $x^*$ in \eqref{equ:3:12} by $\hat x$. Thus the sequence $\{V(x_k,\hat x)\}$ is a convergent sequence in the sense of almost surely. Since $\hat x$ is an accumulation point of $\{x_k\}$, it easily follows that $\{V(x_k,\hat x)\}$ converges to zero, i.e., $\{x_k\}$ converges to $\hat x\in X^*.$ Therefore, the previous conclusion together with \eqref{equ:3:13} implies the convergence of $\{x_{k+1/2}\}.$ This completes the proof.$\hfill\square$

\begin{remark}
This result is notable from several standpoints. (i) Theorem \ref{the:3.1} establishes the a.s. convergence of Algorithm 1 through a more concise and effective way than using the supermartingale convergence theorem, in which we do not require a strict recursive relation or control the upper bound of oracle errors that were considered in \cite{I2019,KS}. (ii) In contrast with \cite{KS}, Theorem \ref{the:3.1} is not necessarily to require the pseudomonotonicity-plus and the information of the Lipschitz constant of the mapping $F$. (iii) The small stepsize policy considered in \cite{KS} is replaced by the adaptive stepsize based on the line search rule, which may bring faster convergence speed. (iv) We do not assume the set $X$ or the mapping $F$ to be bounded, and so this can relax some restrictions of \cite{KS} from the theoretical perspective.
\end{remark}

\section{Convergence rate and complexity}\label{sec:Convergence rate and complexity}
\setcounter{equation}{0}
\noindent
In this section, we investigate the convergence rate and the complexity of Algorithm 1 under two different cases. When $X$ is bounded, the analysis of the convergence rate and the complexity in terms of the gap function are established. In the case that $X$ is unbounded, we discuss the similar properties according to the natural residual function.

Below, we first establish the convergence rate and the complexity of the gap function when $X$ is bounded. Let us define the following notations:
\begin{align*}
C_1(\beta,c):=&\frac{B^2_*}{4}(\frac{\alpha\beta\mathcal{Q}}{\alpha-1}+c\gamma^2_0);\\
\triangle U_k:=&\Big(\frac{\alpha\beta\mathcal{Q}}{\alpha-1}+4\beta\Big)\|\epsilon^k_1\|^2_*+(c\gamma^2_0+4\beta)\|\epsilon^k_2\|^2_*; \\ \triangle V_k:=&\frac{2\alpha\beta\mathcal{Q}}{\alpha-1}F(x_k)^T\epsilon^k_1+2c\gamma^2_kF(x_{k+1/2})^T\epsilon^k_2-\gamma_k(x_k-x^*)^T\epsilon^k_2.
\end{align*}

The following lemma describes a recursive relation which involves two martingale difference sequences.

\begin{lemma}\label{lemma:4.1}
Let Assumptions \ref{ass:2.1} and \ref{ass:3.1} (i) hold. Assume that $X$ is a bounded set and $x^*\in X^*$. Then the sequence $\{x_k\}$ generated by Algorithm 1 satisfies
\begin{align}\label{equ:4:1}
V(x_{k+1},x^*)\leq&\Big(1+\frac{2\gamma^2_k}{\alpha\beta}\Big)V(x_k,x^*)-(\frac{\alpha}{2}-\frac{1}{c})\|x_{k+1}-x_k\|^2-u_k+C_1(\beta,c)+\triangle U_k+\triangle V_k,
\end{align}
where $u_k=\gamma_kF(x_k)^T(x_k-x^*)$ and $c$ is a positive constant.
\end{lemma}
\noindent{\bf Proof.} By \eqref{equ:3:3} and Lemma \ref{lemma:2.3},
\begin{equation*}
V(x_k,x)\geq\gamma_k\widehat F(x_{k+1/2};\xi_{k+1/2})^T(x_{k+1}-x)+V(x_k,x_{k+1})+V(x_{k+1},x),~~\forall x\in X.
\end{equation*}
Adding and subtracting $x_k$ from the first term on the right hand side in above inequality, we have
\begin{align*}%\label{equ:4:3}
V(x_k,x)\geq&\gamma_k\widehat F(x_{k+1/2};\xi_{k+1/2})^T(x_{k}-x)+\gamma_k\widehat F(x_{k+1/2};\xi_{k+1/2})^T(x_{k+1}-x_k)\nonumber\\
&+V(x_k,x_{k+1})+V(x_{k+1},x).
\end{align*}
Letting $x=x^*,$ adding and subtracting $\gamma_kF(x_k)^T(x_k-x^*)$ from the first term on the right in above inequality, we have
\begin{align*}
V(x_k,x^*)\geq&\gamma_kF(x_k)^T(x_k-x^*)+\gamma_k(x_k-x^*)^T\epsilon^k_2+\gamma_k(F(x_{k+1/2})-F(x_k))^T(x_k-x^*)\\
&+\gamma_k\widehat F(x_{k+1/2};\xi_{k+1/2})^T(x_{k+1}-x_k)+V(x_k,x_{k+1})+V(x_{k+1},x^*).
\end{align*}
Rearranging and completing squares, for any positive $\beta$ and $c$, we get
\begin{align}\label{equ:4:4}
V(&x_k,x^*)+\frac{\gamma^2_k}{\beta}\|x_k-x^*\|^2+\beta\|F(x_{k+1/2})-F(x_k)\|^2_*+c\gamma^2_k\|\widehat F(x_{k+1/2};\xi_{k+1/2})\|^2_*+\frac{1}{c}\|x_{k+1}-x_k\|^2\nonumber\\
&\geq u_k+\gamma_k(x_k-x^*)^T\epsilon^k_2+V(x_k,x_{k+1})+V(x_{k+1},x^*).
\end{align}
On the other hand, is follows from \eqref{equ:2:3} and \eqref{equ:3:1} that
\begin{align}\label{equ:4:5}
\|F(x_{k+1/2})-F(x_k)\|^2_*=&\|\widehat{F}(x_{k+1/2};\xi_{k+1/2})-\epsilon^k_2-\widehat{F}(x_{k};\xi_{k})+\epsilon^k_1\|^2_*\nonumber\\
\leq&2\|\widehat{F}(x_{k+1/2};\xi_{k+1/2})-\widehat{F}(x_{k};\xi_{k})\|^2_*+2\|\epsilon^k_1-\epsilon^k_2\|^2_*\nonumber\\
\leq&\frac{2\alpha}{\gamma^2_k}V(x_k,x_{k+1/2})+4\|\epsilon^k_1\|^2_*+4\|\epsilon^k_2\|^2_*\nonumber\\
\leq&\frac{\alpha \mathcal{Q}}{\gamma^2_k}\|x_{k+1/2}-x_k\|_*^2+4\|\epsilon^k_1\|^2_*+4\|\epsilon^k_2\|^2_*.
\end{align}
By \eqref{equ:4:4}, \eqref{equ:4:5} and \eqref{equ:2:2},
\begin{align}\label{equ:4:6}
V(x_{k+1},x^*)\leq&\Big(1+\frac{2\gamma^2_k}{\alpha\beta}\Big)V(x_k,x^*)-(\frac{\alpha}{2}-\frac{1}{c})\|x_{k+1}-x_k\|^2-\gamma_k(x_k-x^*)^T\epsilon^k_2-u_k\nonumber\\
&+\beta\Big(\frac{\alpha\mathcal{Q}}{\gamma^2_k}\|x_{k+1/2}-x_k\|^2_*+4\|\epsilon^k_1\|^2_*+4\|\epsilon^k_2\|^2_* \Big)\nonumber\\
&+c\gamma^2_k\Big(\|F(x_{k+1/2})\|^2_*+\|\epsilon^k_2\|^2_*+2F(x_{k+1/2})^T\epsilon^k_2\Big).
\end{align}
Similarly, by \eqref{equ:3:2} and Lemma \ref{lemma:2.3}, $$\gamma_k(F(x_k)+\epsilon^k_1)^Tx_{k+1/2}+V(x_k,x_{k+1/2})+V(\tilde{x},x_{k+1/2})\leq\gamma_k(F(x_k)+\epsilon^k_1)^Tx_k+V(\tilde{x},x_k), ~~\forall \tilde{x}\in X.$$
Setting $\tilde{x}=x_k$ in the above inequality and using the fact that $V(x_k,x_k)=0$ and $V(x_k,x_{k+1/2})\geq0$, we have
\begin{align*}
\gamma_k(F(x_k)+\epsilon^k_1)^Tx_{k+1/2}+V(x_k,x_{k+1/2})\leq&\gamma_k(F(x_k)+\epsilon^k_1)^Tx_k+V(x_k,x_k)\\
=&\gamma_k(F(x_k)+\epsilon^k_1)^Tx_k.
\end{align*}
It follows from \eqref{equ:2:2} that
\begin{align*}
\frac{\alpha}{2}\|x_{k+1/2}-x_k\|^2\leq& V(x_k,x_{k+1/2})\leq\gamma_k(x_k-x_{k+1/2})^T(F(x_k)+\epsilon^k_1)\\
\leq&\gamma_k\|x_k-x_{k+1/2}\|\cdot\|F(x_k)+\epsilon^k_1\|_*\\
\leq&\frac{1}{2}\|x_k-x_{k+1/2}\|^2+\frac{\gamma^2_k}{2}\|F(x_k)+\epsilon^k_1\|^2_*,
\end{align*}
which together with Assumption \ref{ass:3.1} (i) yields
\begin{align*}
\|x_{k+1/2}-x_k\|^2\leq&\gamma^2_k\frac{\|F(x_k)+\epsilon^k_1\|^2_*}{\alpha-1}\\
=&\frac{\gamma^2_k}{\alpha-1}\Big(\|F(x_k)\|^2_*+\|\epsilon^k_1\|^2_*+2F(x_k)^T\epsilon^k_1\Big)\\
\leq&\frac{\gamma^2_k}{\alpha-1}\Big(\frac{B^2_*}{4}+\|\epsilon^k_1\|^2_*+2F(x_k)^T\epsilon^k_1\Big).
\end{align*}
From \eqref{equ:4:6}, we get
\begin{align*}
V(x_{k+1},x^*)\leq&\Big(1+\frac{2\gamma^2_k}{\alpha\beta}\Big)V(x_k,x^*)-(\frac{\alpha}{2}-\frac{1}{c})\|x_{k+1}-x_k\|^2-\gamma_k(x_k-x^*)^T\epsilon^k_2-u_k\\
&+\frac{2\alpha\beta\mathcal{Q}}{\alpha-1}F(x_k)^T\epsilon^k_1+2c\gamma^2_kF(x_{k+1/2})^T\epsilon^k_2+\frac{B^2_*}{4}(\frac{\alpha\beta\mathcal{Q}}{\alpha-1}+c\gamma^2_0)\\
&+\Big(\frac{\alpha\beta\mathcal{Q}}{\alpha-1}+4\beta\Big)\|\epsilon^k_1\|^2_*+(c\gamma^2_0+4\beta)\|\epsilon^k_2\|^2_*.
\end{align*}
Using the definitions of $\triangle U_k, \triangle V_k$ and $C_1(\beta,c)$, we can get the conclusion. This completes the proof.$\hfill\square$

Next, we aim at computing the conditional expectation value of the terms $\triangle U_k$ and $\triangle V_k$ in the right side of \eqref{equ:4:1}. Since $\{\epsilon^k_1\}$ and $\{\epsilon^k_2\}$ are both martingale difference sequences, we can control them by a relatively straightforward way.

\begin{lemma} \label{lemma:4.2}
Let Assumptions \ref{ass:2.1}, \ref{ass:3.1} (ii) and \ref{ass:3.3} hold. Assume that $X$ is  a bounded set and $x^*\in X^*$. Then the sequence $\{x_k\}$ generated by Algorithm 1 satisfies
$$\mathbb{E}[\triangle V_k|\mathcal{F}_k]=0$$
and
$$\mathbb{E}[\triangle U_k|\mathcal{F}_k]=\frac{2(\frac{\alpha\beta\mathcal{Q}}{\alpha-1}+4\beta)}{N_k}\delta^2(x^*)(1+\|x_{k}-x^*\|^2)
+\frac{(c\gamma^2_0+4\beta)v^2}{N_k},$$
where $\delta(x^*):=\max\{|\|f(x^*,\xi)-F(x^*)\|_*|_2,\mathcal{L}+|L(\xi)|_2\}.$
\end{lemma}

\noindent{\bf Proof.} From the definition of $\{\epsilon_1^k\}$ and the facts that $x_k\in \mathcal{F}_k, \;\xi_{j,k} \perp\perp\mathcal{F}_k$, we get
\begin{align}\label{equ:4:7}
\mathbb{E}[\|\epsilon^k_1\|^2|\mathcal{F}_k]=&\mathbb{E}[\|\frac{1}{N_k}\sum^{N_k}_{j=1}f(x_k,\xi_{j,k})-F(x_k)\|^2|\mathcal{F}_k]\nonumber\\
=&V\Big(\frac{1}{N_k}\sum^{N_k}_{j=1}f(x_k,\xi_{j,k})|\mathcal{F}_k\Big)=\frac{1}{N_k^2}V\Big(\sum^{N_k}_{j=1}f(x_k,\xi_{j,k})|\mathcal{F}_k\Big)\nonumber\\
=&\frac{1}{N_k^2}\sum^{N_k}_{j=1}V(f(x_k,\xi_{j,k})|\mathcal{F}_k)=\frac{1}{N_k^2}\sum^{N_k}_{j=1}V(f(x_k,\xi)|\mathcal{F}_k)\nonumber\\
=&\frac{1}{N_k}V(f(x_k,\xi)|\mathcal{F}_k)=\frac{1}{N_k}|\|f(x_k,\xi)-F(x_k)\|_*|^2_2,
\end{align}
where the third equality is based on the property of the variance, the fourth equality comes from the independence of $\xi_{j,k}$ and the fifth equality is due to the fact that all $\xi_{j,k}$ obey the same distribution. By Assumption 3.8 of \cite{I2017} and the definition of $\delta(x^*)$,
\begin{align*}
|\|f(x_k,\xi)-F(x_k)\||_2&\leq|\|f(x_k,\xi)-f(x^*,\xi)\||_2+|\|f(x^*,\xi)-F(x^*)\|_*|_2+\|F(x_k)-F(x^*)\|\\
&\leq|\|f(x^*,\xi)-F(x^*)\|_*|_2+(\mathcal{L}+|L(\xi)|_2)\|x_k-x^*\|\\
&\leq\delta(x^*)(1+\|x_k-x^*\|),
\end{align*}
which together with \eqref{equ:4:7} yields
\begin{equation}\label{equ:4:8}
\mathbb{E}[\|\epsilon^k_1\|_*^2|\mathcal{F}_k]\leq\frac{1}{N_k}\delta^2(x^*)(1+\|x_k-x^*\|)^2\leq\frac{2}{N_k}\delta^2(x^*)(1+\|x_k-x^*\|^2).
\end{equation}
Similarly, using the definition of $\{\epsilon^k_2\}$ and the facts that $x_{k+1/2}\in \mathcal{F}_{k+1/2}, \;\xi_{j,k+1/2}\perp\perp\mathcal{F}_{k+1/2}$, we can obtain
\begin{align*}
\mathbb{E}[\|\epsilon^k_2\|_*^2|\mathcal{F}_{k+1/2}]
=\frac{1}{N_k}|\|f(x_{k+1/2},\xi)-F(x_{k+1/2})\|_*|^2_2.
\end{align*}
It follows from Assumption \ref{ass:3.1} (ii) that
\begin{equation}\label{equ:4:9}
\mathbb{E}[\|\epsilon^k_2\|_*^2|\mathcal{F}_{k}]=\mathbb{E}[\mathbb{E}[\|\epsilon^k_2\|_*^2|\mathcal{F}_{k+1/2}]|\mathcal{F}_k]=\frac{1}{N_k}\mathbb{E}[v^2|\mathcal{F}_k]=\frac{v^2}{N_k}.
\end{equation}
Combining \eqref{equ:4:8} and \eqref{equ:4:9} yields
\begin{align*}
\mathbb{E}[\triangle U_k|\mathcal{F}_k]&=(\frac{\alpha\beta\mathcal{Q}}{\alpha-1}+4\beta)\mathbb{E}[\|\epsilon^k_1\|^2_*|\mathcal{F}_k]+(c\gamma^2_0+4\beta)\mathbb{E}[\|\epsilon^k_2\|^2_*|\mathcal{F}_k]\\
&=\frac{2(\frac{\alpha\beta\mathcal{Q}}{\alpha-1}+4\beta)}{N_k}\delta^2(x^*)(1+\|x_{k}-x^*\|^2)+\frac{(c\gamma^2_0+4\beta)v^2}{N_k}.
\end{align*}
By \eqref{equ:3:4} and $x_k\in\mathcal{F}_k,x_{k+1/2}\in\mathcal{F}_{k+1/2}$ and $\mathcal{F}_k\subseteq\mathcal{F}_{k+1/2}$, one has
$$\mathbb{E}[\triangle V_k|\mathcal{F}_k]=\frac{2\alpha\beta\mathcal{Q}}{\alpha-1}F(x_k)^T\mathbb{E}[\epsilon^k_1|\mathcal{F}_k]+2c\gamma^2_k\mathbb{E}[F(x_{k+1/2})^T\mathbb{E}[\epsilon^k_2|\mathcal{F}_{k+1/2}]|\mathcal{F}_k]-\gamma_k(x_k-x^*)^T\mathbb{E}[\epsilon^k_2|\mathcal{F}_k]=0.$$
This completes the proof.$\hfill\square$

\begin{remark}
In Lemma \ref{lemma:4.2}, we obtain the upper bounds of $\|\epsilon^k_1\|^2$ and $\|\epsilon^k_2\|^2$ by the uniformly bounded variance assumption and its improved form, which were initially studied in \cite{I2017}. For every iteration $k$, we adopt the multisample strategy to approximate the value of $F$ while the algorithm considered in \cite{KS} only uses a single sample, which may lead to the high iterative complexity and a slow convergence speed in practice.
\end{remark}

It is worth mentioning that $\mathbb{E}[V(x_k,x^*)]$ is bounded when $X$ is bounded, which implies that there exists a constant $M$ satisfying $\sup_{k\geq0}\mathbb{E}[V(x_k,x^*)] \leq M.$ Further, in order to establish the convergence rate of Algorithm 1,
we need the following gap function (see \cite{HP,H} and references therein):
\begin{equation}\label{equ:4:10}
g(x):=\sup\limits_{z\in X}\langle F(x),x-z\rangle.
\end{equation}
Note that $g(x)\geq0$ for any $x\in X$ and the point $x^*\in X$ is a solution of the SVI if and only if $g(x^*)=0.$

To ease notation, we make the following definitions.
$$C_2(\beta,x^*):=\frac{4\Big(\frac{\alpha\beta\mathcal{Q}}{\alpha-1}+4\beta\Big)\delta^2(x^*)}{\alpha}; ~~ C_3(\beta,c,v,x^*):=2\Big(\frac{\alpha\beta\mathcal{Q}}{\alpha-1}+4\beta\Big)\delta^2(x^*)+(c\gamma^2_0+4\beta)v^2.$$

The convergence rate of Algorithm 1 in terms of the gap function is given as follows.

\begin{theorem} \label{the:4.1}
Let Assumptions \ref{ass:2.1}, \ref{ass:3.1} and \ref{ass:3.3} hold. Assume that the sequence $\{x_k\}$ is  generated by Algorithm 1 and $X$ is a bounded set. Let $\beta>0$ and $c>0$ be chosen such that $\frac{\alpha}{2}-\frac{1}{c}\geq0$. Then for any $x^*\in X^*$ and $k\geq0$, the following statements are true:
\begin{itemize}
\item [(i)] If $N_k=\mathcal{N}\lceil (k+\lambda)(\ln (k+\lambda))^{1+b}\rceil$ with $\mathcal{N}=\mathcal{O}(d), \lambda>2$ and $b>0,$ then
\begin{align}\label{equ:4:11}
\min_{i\in\{0,\cdots,k\}}\mathbb{E}[g(x_i)]
\leq& \frac{|L(\xi,\xi_{+1/2})|_2}{(k+1)\min\{\frac{\alpha\theta}{\sqrt{2}},\gamma_0\}}
\Big[V(x_0,x^*)+\frac{2(k+1)M\gamma^2_0}{\alpha\beta}+\frac{MC_2(\beta,x^*)}
{\mathcal{N}b(\ln(\lambda-1))^b}\nonumber\\
&+(k+1)C_1(\beta,c)+\frac{C_3(\beta,c,v,x^*)}{\mathcal{N}b(\ln(\lambda-1))^b}\Big];
\end{align}
\item [(ii)]  If $N_k=\lceil \frac{(k+1)^{1.5}}{d}\rceil$, then
\begin{align*}
\min_{i\in\{0,\cdots,k\}}\mathbb{E}[g( x_i)]
\leq&\frac{|L(\xi,\xi_{+1/2})|_2}{(k+1)\min\{\frac{\alpha\theta}{\sqrt{2}},\gamma_0\}}
\Big[V(x_0,x^*)+\frac{2(k+1)M\gamma^2_0}{\alpha\beta}+\frac{MC_2(\beta,x^*)}{3d}\\
&+(k+1)C_1(\beta,c)+\frac{C_3(\beta,c,v,x^*)}{3d}\Big].
\end{align*}
\end{itemize}
\end{theorem}

\noindent{\bf Proof.} (i) Taking $\mathbb{E}[\cdot|\mathcal{F}_k]$ in \eqref{equ:4:1} and by Lemmas \ref{lemma:3.2} and \ref{lemma:4.2},
\begin{align*}
\mathbb{E}[V(x_{k+1},x^*)|\mathcal{F}_k]\leq&(1+\frac{2\gamma^2_k}{\alpha\beta})V(x_k,x^*)-\mathbb{E}[\gamma_k|\mathcal{F}_k]F(x_k)^T(x_k-x^*)\\
&+C_1(\beta,c)+\mathbb{E}[\triangle U_k|\mathcal{F}_k]+\mathbb{E}[\triangle V_k|\mathcal{F}_k]\\
\leq&(1+\frac{2\gamma^2_0}{\alpha\beta})V(x_k,x^*)-\frac{\min \{\frac{\alpha\theta}{\sqrt{2}},\gamma_0\}}{|L(\xi,\xi_{+1/2})|_2}F(x_k)^T(x_k-x^*)\\
&+C_1(\beta,c)+\frac{2(\frac{\alpha\beta\mathcal{Q}}{\alpha-1}
+4\beta)}{N_k}\delta^2(x^*)(1+\|x_k-x^*\|^2)+\frac{(c\gamma^2_0+4\beta)v^2}{N_k}.
\end{align*}
This together with \eqref{equ:2:2} yields
\begin{align*}
\mathbb{E}[V(x_{k+1},x^*)|\mathcal{F}_k]
\leq&(1+\frac{2\gamma^2_0}{\alpha\beta})V(x_k,x^*)-\frac{\min \{\frac{\alpha\theta}{\sqrt{2}},\gamma_0\}}{|L(\xi,\xi_{+1/2})|_2}F(x_k)^T(x_k-x^*)\\
&+C_1(\beta,c)+\frac{2(\frac{\alpha\beta\mathcal{Q}}{\alpha-1}
+4\beta)}{N_k}\delta^2(x^*)\Big(1+\frac{2}{\alpha}V(x_k,x^*)\Big)+\frac{(c\gamma^2_0+4\beta)v^2}{N_k}\\
\leq&(1+\frac{2\gamma^2_0}{\alpha\beta})V(x_k,x^*)-\frac{\min \{\frac{\alpha\theta}{\sqrt{2}},\gamma_0\}}{|L(\xi,\xi_{+1/2})|_2}F(x_k)^T(x_k-x^*)\\
&+C_1(\beta,c)+C_2(\beta,x^*)\frac{V(x_k,x^*)}{N_k}+\frac{C_3(\beta,c,v,x^*)}{N_k}.
\end{align*}
From \eqref{equ:4:10}, we obtain
\begin{align*}%\label{equ:4:12}
\frac{\min\{\frac{\alpha\theta}{\sqrt{2}},\gamma_0\}}{|L(\xi,\xi_{+1/2})|_2}g(x_k)
\leq&(1+\frac{2\gamma^2_0}{\alpha\beta})V(x_k,x^*)-\mathbb{E}[V(x_{k+1},x^*)|\mathcal{F}_k]\nonumber\\
&+C_1(\beta,c)+C_2(\beta,x^*)\frac{V(x_k,x^*)}{N_k}+\frac{C_3(\beta,c,v,x^*)}{N_k}.
\end{align*}
Taking expectation in above inequality, using $\mathbb{E}[\mathbb{E}[\cdot|\mathcal{F}_i]]=\mathbb{E}[\cdot]$ and summing it from $i=0$ to $i=k$, we have
\begin{align*}
\frac{\min\{\frac{\alpha\theta}{\sqrt{2}},\gamma_0\}}{|L(\xi,\xi_{+1/2})|_2}\sum^{k}_{i=0}\mathbb{E}[g(x_i)]\leq& V(x_0,x^*)+\frac{2\gamma^2_0}{\alpha\beta}\sum^{k}_{i=0}\mathbb{E}[V(x_i,x^*)]+(k+1)C_1(\beta,c)\\
&+C_2(\beta,x^*)\sum^{k}_{i=0}\frac{\mathbb{E}[V(x_i,x^*)]}{N_i}+C_3(\beta,c,v,x^*)\sum^{k}_{i=0}\frac{1}{N_i}\\
\leq&V(x_0,x^*)+\frac{2M(k+1)\gamma^2_0}{\alpha\beta}+(k+1)C_1(\beta,c)\\
&+C_2(\beta,x^*)\sum^{k}_{i=0}\frac{M}{N_i}+C_3(\beta,c,v,x^*)\sum^{k}_{i=0}\frac{1}{N_i}.
\end{align*}
This implies that
\begin{align*}
\sum^{k}_{i=0}\mathbb{E}[g(x_i)]\leq&\frac{|L(\xi,\xi_{+1/2})|_2}{\min\{\frac{\alpha\theta}{\sqrt{2}},\gamma_0\}}
\Big[V(x_0,x^*)+\frac{2M(k+1)\gamma^2_0}{\alpha\beta}+(k+1)C_1(\beta,c)\\
&+C_2(\beta,x^*)\sum^{k}_{i=0}\frac{M}{N_i}+C_3(\beta,c,v,x^*)\sum^{k}_{i=0}\frac{1}{N_i}\Big].
\end{align*}
By dividing both sides of the above inequality by $k+1$,
\begin{align}\label{equ:4:13}
\min_{i\in\{0,\cdots,k\}}\mathbb{E}[g(x_i)]\leq&\frac{1}{k+1}\sum^{k}_{i=0}
\mathbb{E}[g(x_i)]\nonumber\\
\leq& \frac{|L(\xi,\xi_{+1/2})|_2}{(k+1)\min\{\frac{\alpha\theta}{\sqrt{2}},\gamma_0\}}
\Big[V(x_0,x^*)+\frac{2M(k+1)\gamma^2_0}{\alpha\beta}+(k+1)C_1(\beta,c)\nonumber\\
&+C_2(\beta,x^*)\sum^{k}_{i=0}\frac{M}{N_i}+C_3(\beta,c,v,x^*)\sum^{k}_{i=0}
\frac{1}{N_i}\Big].
\end{align}
From the definition of $N_k$, we get
\begin{align}\label{equ:4:14}
\sum^k_{i=0}\frac{1}{N_i}\leq\sum^{\infty}_{i=0}\frac{1}{N_i}\leq\int^{\infty}_{-1}\frac{d_q}{\mathcal{N}(q+\lambda)(\ln(q+\lambda))^{1+b}}=\frac{1}{\mathcal{N}b(\ln(\lambda-1))^b}.
\end{align}
Combining \eqref{equ:4:13} and \eqref{equ:4:14} yields the conclusion.

(ii) Letting $N_k=\lceil \frac{(k+1)^{1.5}}{d}\rceil$, we have
\begin{align}\label{equ:4:15}
\sum^k_{i=0}\frac{1}{N_i}\leq\frac{1}{N_0}+\sum^k_{i=0}\frac{1}{N_i}\leq d+\sum^k_{i=0}\frac{1}{N_i}\leq d+d\int^{\infty}_{0}\frac{d_q}{(q+1)^{1.5}}=3d.
\end{align}
By \eqref{equ:4:13}, we can get the conclusion. The proof is completed.$\hfill\square$

Based on Theorem \ref{the:4.1}, we have the following results about the iteration complexity and the oracle complexity of Algorithm 1.

\begin{cor} \label{cor:4.1}
Let Assumptions \ref{ass:2.1}, \ref{ass:3.1} and \ref{ass:3.3} hold. Assume that $X$ is a bounded set and $x^*\in X^*$. Set $N_k=\mathcal{N}\lceil (k+\lambda)(\ln (k+\lambda))^{1+b}\rceil$ or
$N_k=\lceil \frac{(k+1)^{1.5}}{d}\rceil$. Then the following statements are true:
\begin{itemize}
\item [(i)] In terms of $\min\limits_{i\in\{0,\cdots,k\}}\mathbb{E}[g(x_i)]$, the iteration complexity of Algorithm 1 is $\mathcal{O}(\frac{1}{\varepsilon}$);

\item [(ii)] Both in the sense of almost surely and in the sense of expectation, the oracle complexity of Algorithm 1 is $\mathcal{O}(\frac{1}{\varepsilon^2})$.
\end{itemize}
\end{cor}

\noindent{\bf Proof.} Let $N_k=\mathcal{N}\lceil (k+\lambda)(\ln (k+\lambda))^{1+b}\rceil$.

(i) If $$\frac{|L(\xi,\xi_{+1/2})|_2}{(k+1)\min\{\frac{\alpha\theta}{\sqrt{2}},\gamma_0\}}
\Big[V(x_0,x^*)+\frac{2(k+1)M\gamma^2_0}{\alpha\beta}+\frac{MC_2(\beta,x^*)}{\mathcal{N}b(\ln(\lambda-1))^b}
+(k+1)C_1(\beta,c)+\frac{C_3(\beta,c,v,x^*)}{\mathcal{N}b(\ln(\lambda-1))^b}\Big]\leq \varepsilon$$
or
$$k+1\geq \frac{|L(\xi,\xi_{+1/2})|_2}{\varepsilon\min\{\frac{\alpha\theta}{\sqrt{2}},\gamma_0\}}
\Big[V(x_0,x^*)+\frac{2(k+1)M\gamma^2_0}{\alpha\beta}+\frac{MC_2(\beta,x^*)}{\mathcal{N}b(\ln(\lambda-1))^b}
+(k+1)C_1(\beta,c)+\frac{C_3(\beta,c,v,x^*)}{\mathcal{N}b(\ln(\lambda-1))^b}\Big]$$
then by \eqref{equ:4:11},
$$\min\limits_{i\in\{0,\cdots,k\}}\mathbb{E}[g(x_i)]\leq \varepsilon.$$
This implies that the iteration complexity is $\mathcal{O}(\frac{1}{\varepsilon})$.

(ii) Let $l_i$ be the number of the line search in the iteration $i$. After $k$ iterations, the oracle complexity is upper bounded by
\begin{align}\label{equ:4:16}
\sum^k_{i=0}(1+l_i)N_i&\leq\{1+\max_{i\in\{0,\cdots,k\}}l_i\}\sum^k_{i=0}\Big(\mathcal{N}(i+\lambda)(\ln (i+\lambda))^{1+b}+1\Big)\nonumber\\
&=\{1+\max_{i\in\{0,\cdots,k\}}l_i\}\sum^{k+\lambda}_{i=\lambda}\Big(\mathcal{N}i(\ln (i))^{1+b}+1\Big)\nonumber\\
&\leq \{1+\max_{i\in\{0,\cdots,k\}}l_i\}\sum^{k+\lambda}_{i=\lambda}\Big(\mathcal{N}i(\ln (k+\lambda))^{1+b}+1\Big)\nonumber\\
&=\{1+\max_{i\in\{0,\cdots,k\}}l_i\}\Big(\mathcal{N}(\ln (k+\lambda))^{1+b}\frac{(k+1)(2\lambda+k)}{2}+k+1\Big).
\end{align}
It follows from Lemma \ref{lemma:3.2} that for all $k$,
$$\min\{\frac{\alpha\theta}{\bar{L}(\xi_k,\xi_{k+1/2})\sqrt{2}}, \gamma_0\}\leq \gamma_k\leq\gamma_0\theta^{l_k},\;\;a.s.$$
Then
$$l_k\leq \log_{\frac{1}{\theta}}\Big(\frac{\gamma_0}{\min\{\frac{\alpha\theta}{\bar {L}(\xi_k,\xi_{k+1/2})\sqrt{2}},\gamma_0\}}\Big)\leq \log_{\frac{1}{\theta}}\Big(\frac{\gamma_0 \bar {L}(\xi_k,\xi_{k+1/2})}{\min\{\frac{\alpha\theta}{\sqrt{2}},\gamma_0\}}\Big),\;\;a.s.$$
From \eqref{equ:4:16} and the iteration complexity $\mathcal{O}(\frac{1}{\varepsilon})$, we have a.s. the oracle complexity is $\mathcal{O}(\frac{1}{\epsilon^2})$. In view of the concavity of $\log_{\frac{1}{\theta}}t$ and Lemma \ref{lemma:2.2}, we
obtain
\begin{equation}\label{equ:4:17}
\mathbb{E}[l_k]\leq\mathbb{E}\Big[\log_{\frac{1}{\theta}}\Big(\frac{\gamma_0 \bar {L}(\xi_k,\xi_{k+1/2})}{\min\{\frac{\alpha\theta}{\sqrt{2}},\gamma_0\}}\Big)\Big]\leq\log_{\frac{1}{\theta}}\Big(\frac{\gamma_0 \mathcal{L}}{\min\{\frac{\alpha\theta}{\sqrt{2}},\gamma_0\}}\Big).
\end{equation}
Taking expectation in \eqref{equ:4:16} and using \eqref{equ:4:17}, we have that in the sense of expectation, the oracle complexity is $\mathcal{O}(\frac{1}{\varepsilon^2}).$

(ii) The analysis for $N_k=\lceil \frac{(k+1)^{1.5}}{d}\rceil$ is similar, so we omit it.
$\hfill\square$

\begin{remark}
\begin{itemize}
\item [(i)] In \cite{KS}, Kannan and Shanbhag derived the same convergence rate when the mapping $F$ is strongly pseudomonotone. However, Theorem \ref{the:4.1} does not require this assumption.

\item [(ii)] Corollary \ref{cor:4.1} shows the complexity results of Algorithm 1. Although there is no complexity analysis of the algorithm in \cite{KS}, we can easily get that its iteration and oracle complexities of the merit function are both $\mathcal{O}(1/\varepsilon^2)$. Hence, our proposed algorithm may have better iteration complexity than the one of \cite{KS}.
\end{itemize}
\end{remark}

Given the expensive costs of function calculations, a potential improvement is to reduce the oracle calls as reasonably as possible while retaining the algorithm's properties. In what follows, we will discuss the convergence rate and complexity of the natural residual function when $X$ is unbounded. Note also that in this framework, only a single oracle call is required at each iteration.

We start with some pivotal lemmas.

\begin{lemma} \label{lemma:4.3}  \cite{NA}
Given $x\in X$ and $\gamma>0$, then for all $r_1, r_2\in\mathbb{R}^n$, we have
$$\|\mathbf{P}_X(x,r_1)-\mathbf{P}_X(x,r_1)\|\leq\frac{\gamma}{\alpha}\|r_1-r_2\|_*.$$
\end{lemma}

\begin{lemma} \label{lemma:4.4}
Given $x\in X$ and $r\in\mathbb{R}^n$, denote $x^+:=\mathbf{P}_X(x,r)$. Then, for any $u\in X$, we have
$$\big(\nabla{s(x^+)}-\nabla{s(x)}\big)^T(x^+-u)=V(x^+,u)+V(x,x^+)-V(x,u).$$
\end{lemma}

\noindent{\bf Proof.} From \eqref{equ:2:1}, we have
\begin{align*}
V(x^+,u)-V(x,u)=&s(u)-s(x^+)-\nabla{s(x^+)}^T(u-x^+)-\big[s(u)-s(x)-\nabla{s(x)}^T(u-x)\big]\\
=&s(x)-s(x^+)+\nabla{s(x^+)}^T(x^+-u)-\nabla{s(x)}^T(x^+-u)+\nabla{s(x)}^T(x^+-x)\\
=&\big(\nabla{s(x^+)}-\nabla{s(x)}\big)^T(x^+-u)+s(x)-s(x^+)+\nabla{s(x)}^T(x^+-x)\\
=&\big(\nabla{s(x^+)}-\nabla{s(x)}\big)^T(x^+-u)-V(x,x^+).
\end{align*}
This complete the proof. $\hfill\square$

To analyze the convergence rate of Algorithm 1, we give an important recursive relation in the sense of the expectation. To ease notation, we make the following definitions:
$$C_4(x^*):=\frac{\gamma_0^2\delta^2(x^*)}{2\alpha},~~ C_5(x^*):=\frac{\gamma_0^2\delta^2(x^*)}{\alpha^2}.$$

\begin{lemma} \label{lemma:4.5}
Let Assumptions \ref{ass:2.1}, \ref{ass:3.2} and \ref{ass:3.3} hold. Then for any $k$ and $x^*\in X^*$, the sequence $\{x_k\}$ generated by Algorithm 1 satisfies
\begin{align}\label{equ:4:18}
\mathbb{E}[V(x_{k+1},x^*)|\mathcal{F}_k]\leq&\mathbb{E}[V(x_{k},x^*)|\mathcal{F}_k]+\frac{C_4(x^*)}{N_k}\Big(1+\frac{2}{\alpha}V(x_k,x^*)\Big)-\frac{\alpha}{8}R^2_{\gamma_k}(x_k,F(x_k)).
\end{align}
\end{lemma}
\noindent{\bf Proof.} By Lemma \ref{lemma:2.3} and \eqref{equ:3:3},
\begin{align}\label{equ:4:19}
V(x_{k+1},x^*)\leq&V(x_k,x^*)-\gamma_k\widehat F(x_{k+1/2};\xi_{k+1/2})^T(x_{k+1}-x^*)-V(x_k,x_{k+1})\nonumber\\
\leq&V(x_k,x^*)-\gamma_k\widehat F(x_{k+1/2};\xi_{k+1/2})^T(x_{k+1}-x_{k+1/2})\nonumber\\
&-\gamma_k\widehat F(x_{k+1/2};\xi_{k+1/2})^T(x_{k+1/2}-x^*)-V(x_k,x_{k+1})\nonumber\\
=&V(x_k,x^*)-\gamma_k\big(\widehat F(x_{k+1/2};\xi_{k+1/2})-\widehat F(x_k;\xi_k)\big)^T(x_{k+1}-x_{k+1/2})\nonumber\\
&-\gamma_k\widehat F(x_k;\xi_k)^T(x_{k+1}-x_{k+1/2})-\gamma_k\widehat F(x_{k+1/2};\xi_{k+1/2})^T(x_{k+1/2}-x^*)\nonumber\\
&-V(x_k,x_{k+1}).
\end{align}
Using \eqref{equ:3:2} and the optimality condition of \eqref{equ:2:4}, we get
\begin{equation}\label{equ:4:20}
\gamma_k\widehat F(x_k;\xi_k)^T(x_{k+1/2}-x_{k+1})\leq\big(\nabla{s(x_{k+1/2})}-\nabla{s(x_k)}\big)^T(x_{k+1}-x_{k+1/2}).
\end{equation}
By Lemma \ref{lemma:4.4} and \eqref{equ:3:2},
$$\big(\nabla{s(x_{k+1/2})}-\nabla{s(x_k)}\big)^T(x_{k+1/2}-x_{k+1})=V(x_{k+1/2},x_{k+1})+V(x_k,x_{k+1/2})-V(x_k,x_{k+1}).$$
This together with \eqref{equ:4:19} and \eqref{equ:4:20} yields
\begin{align*}%\label{equ:4:21}
V(x_{k+1},x^*)\leq&V(x_k,x^*)-\gamma_k\big(\widehat F(x_{k+1/2};\xi_{k+1/2})-\widehat F(x_k;\xi_k)\big)^T(x_{k+1}-x_{k+1/2})\nonumber\\&+V(x_k,x_{k+1})-V(x_{k+1/2},x_{k+1})-V(x_k,x_{k+1/2})\nonumber\\
&-\gamma_k\widehat F(x_{k+1/2};\xi_{k+1/2})^T(x_{k+1/2}-x^*)-V(x_k,x_{k+1})\nonumber\\
=&V(x_k,x^*)-\gamma_k\widehat F(x_{k+1/2};\xi_{k+1/2})^T(x_{k+1/2}-x^*)-V(x_{k+1/2},x_{k+1})\nonumber\\
&-V(x_k,x_{k+1/2})-\gamma_k\big(\widehat F(x_{k+1/2};\xi_{k+1/2})-\widehat F(x_k;\xi_k)\big)^T(x_{k+1}-x_{k+1/2}).
\end{align*}
On the other hand,
\begin{align*}
-\gamma_k\big(&\widehat F(x_{k+1/2};\xi_{k+1/2})-\widehat F(x_k;\xi_k)\big)^T(x_{k+1}-x_{k+1/2})\\
\leq&\gamma_k\|\widehat F(x_{k+1/2};\xi_{k+1/2})-\widehat F(x_k;\xi_k)\|\cdot\|x_{k+1}-x_{k+1/2}\|\\
\leq&\frac{\gamma_k^2}{2\alpha}\|\widehat F(x_{k+1/2};\xi_{k+1/2})-\widehat F(x_k;\xi_k)\|^2+\frac{\alpha}{2}\|x_{k+1}-x_{k+1/2}\|^2\\
\leq&\frac{1}{2}V(x_k,x_{k+1/2})+V(x_{k+1/2},x_{k+1}).
\end{align*}
Combining the above two observations, we have
\begin{align}\label{equ:4:22}
V(x_{k+1},x^*)\leq V(x_k,x^*)-\gamma_k\widehat F(x_{k+1/2};\xi_{k+1/2})^T(x_{k+1/2}-x^*)-\frac{1}{2}V(x_k,x_{k+1/2}).
\end{align}
Using \eqref{equ:2:5} and Lemma \ref{lemma:4.3}, we have
\begin{align*}
R^2(x_k,F(x_k))=&\|x_k-\mathbf{P}_X(x_k,\gamma_kF(x_k))\|^2\\
=&2\|x_k-x_{k+1/2}\|^2+2\|\mathbf{P}_X(x_k,\gamma_k\widehat F(x_k;\xi_k))-\mathbf{P}_X(x_k,\gamma_kF(x_k))\|^2\\
\leq&2\|x_k-x_{k+1/2}\|^2+\frac{2\gamma_k^2}{\alpha^2}\|\varepsilon^k_1\|^2_*.
\end{align*}
This together with \eqref{equ:2:2} yields
$$-\frac{1}{2}V(x_k,x_{k+1/2})\leq-\frac{\alpha}{4}\|x_k-x_{k+1/2}\|^2\leq\frac{\gamma_k^2}{4\alpha}\|\varepsilon^k_1\|^2_*-\frac{\alpha}{8}R^2_{\gamma_k}(x_k,F(x_k)).$$
Adding it into \eqref{equ:4:22}, we get
\begin{align*}
V(x_{k+1},x^*)\leq&V(x_k,x^*)-\gamma_k\widehat F(x_{k+1/2};\xi_{k+1/2})^T(x_{k+1/2}-x^*)+\frac{\gamma_k^2}{4\alpha}\|\varepsilon^k_1\|^2_*-\frac{\alpha}{8}R^2_{\gamma_k}(x_k,F(x_k))\\
=&V(x_k,x^*)-\gamma_k\big(F(x_{k+1/2})+\varepsilon^k_2\big)^T(x_{k+1/2}-x^*)+\frac{\gamma_k^2}{4\alpha}\|\varepsilon^k_1\|^2_*-\frac{\alpha}{8}R^2_{\gamma_k}(x_k,F(x_k)).
\end{align*}
From Assumption \ref{ass:3.2}, we have that $F(x_{k+1/2})^T(x_{k+1/2}-x^*)\geq0$. It follows that
$$V(x_{k+1},x^*)\leq V(x_k,x^*)-\gamma_k(x_{k+1/2}-x^*)^T\varepsilon^k_2+\frac{\gamma_k^2}{4\alpha}\|\varepsilon^k_1\|^2_*-\frac{\alpha}{8}R^2_{\gamma_k}(x_k,F(x_k)).$$
Taking $\mathbb{E}[\cdot|\mathcal{F}_k]$ in above inequality, and from \eqref{equ:3:4} and \eqref{equ:4:8}, we have
\begin{align*}
\mathbb{E}[V(x_{k+1},x^*)|\mathcal{F}_k]\leq&\mathbb{E}[V(x_{k},x^*)|\mathcal{F}_k]+\frac{C_4(x^*)}{N_k}\Big(1+\frac{2}{\alpha}V(x_k,x^*)\Big)-\frac{\alpha}{8}R^2_{\gamma_k}(x_k,F(x_k)).
\end{align*}
This complete the proof. $\hfill\square$

\begin{remark}
In Lemma \ref{lemma:4.5}, we obtain a desired recursive relation in which it demands to compute the upper bound of one oracle error, while the algorithms considered in \cite{I2019, KS} require the evaluations of two or three oracle errors at each iteration. Based on this, Algorithm 1 may save more calculating time in practice.
\end{remark}

Clearly, the boundedness of $\mathbb{E}[V(x_k,x^*)]$ plays a key role in the convergence rate and complexity of the algorithm. If $X$ is bounded, then we can easily get the conclusion that $\mathbb{E}[V(x_k,x^*)]$ is bounded. However, when $X$ is unbounded, we have the following lemma.

\begin{lemma}\label{lemma:4.6}
Let Assumptions \ref{ass:2.1}, \ref{ass:3.2} and \ref{ass:3.3} hold. Assume that the sequence $\{x_k\}$ is generated by Algorithm 1 and $x^*\in X^*$. Then $\mathbb{E}[V(x_k,x^*)]$ is bounded, i.e., there exists a constant $M$ such that $\sup_{k\geq0}\mathbb{E}[V(x_k, x^*)] \leq M.$
\end{lemma}

\noindent{\bf Proof.}
Taking $\phi\in(0,1)$ and choosing $k_0:= k_0(\phi,C_5(x^*))\in \mathbb{N}$ such that
\begin{align}\label{equ:4:23}
\sum^{\infty}_{i=k_0}\frac{1}{N_i}\leq\frac{\phi}{C_5(x^*)}.
\end{align}
Note that $k_0$ always exists due to $\sum^{\infty}_{k=0}\frac{1}{N_k}<\infty$.

Taking expectation in \eqref{equ:4:18}, using $\mathbb{E}[\mathbb{E}[\cdot]|\mathcal{F}]]=\mathbb{E}[\cdot]$ and summing it from $i=k_0$ to $i=k-1$, we get
\begin{align}\label{equ:4:24}
\mathbb{E}[V(x_{k},x^*)]\leq&\mathbb{E}[V(x_{k_0},x^*)]+C_4(x^*)\sum^{k-1}_{i=k_0}\frac{1}{N_i}+C_5(x^*)\sum^{k-1}_{i=k_0}\frac{\mathbb{E}[V(x_i,x^*)]}{N_i}.
\end{align}
We next prove $\mathbb{E}[V(x_k,x^*)]$ is bounded. In fact, suppose by contradiction that $\mathbb{E}[V(x_k,x^*)]$ is unbounded. Thus for any $a>0$,$$k_a:=\inf \{k\geq k_0:\sqrt{\mathbb{E}[V(x_k,x^*)]}>a\}$$ always exists.
If $k_a=k_0,$ then $a <\sqrt{\mathbb{E}[V(x_k,x^*)]}$, which contradicts to the arbitrariness
of $a$. If $k_a > k_0$, then from \eqref{equ:4:23}, \eqref{equ:4:24} and the definition of $k_a$, we have
\begin{align*}
a^2<&\mathbb{E}[V(x_{k_a},x^*)]\\
\leq&\mathbb{E}[V(x_{k_0},x^*)]+C_4(x^*)\sum^{k-1}_{i=k_0}\frac{1}{N_i}+C_5(x^*)\sum^{k-1}_{i=k_0}\frac{\mathbb{E}[V(x_i,x^*)]}{N_i}\\
\leq&\mathbb{E}[V(x_{k_0},x^*)]+C_4(x^*)\sum^{k-1}_{i=k_0}\frac{1}{N_i}+C_5(x^*)\sum^{k-1}_{i=k_0}\frac{a^2}{N_i}\\
\leq&\mathbb{E}[V(x_{k_0},x^*)]+\frac{C_4(x^*)}{C_5(x^*)}\phi+\phi a^2,
\end{align*}
which together with $\phi\in(0,1)$ yields that
$$a^2<\frac{\mathbb{E}[V(x_{k_0},x^*)]+\frac{C_4(x^*)}{C_5(x^*)}\phi}{1-\phi}.$$
This gives a contradiction with the arbitrariness of $a$. This completes the proof.$\hfill\square$

The convergence rate of Algorithm 1 in terms of the natural residual function is stated as follows.

\begin{theorem} \label{the:4.2}
Let Assumptions \ref{ass:2.1} and \ref{ass:3.1}-\ref{ass:3.3} hold. Assume that the sequence $\{x_k\}$ is generated by Algorithm 1. Then for any $x^*\in X^*$ and $k\geq0$, the following statements are true:
\begin{itemize}
\item [(i)] If $N_k=\mathcal{N}\lceil (k+\lambda)(\ln (k+\lambda))^{1+b}\rceil$, with $\mathcal{N}=\mathcal{O}(d),\lambda>2$ and $b>0,$ then
$$\min_{i\in\{0,\cdots,k\}}\mathbb{E}[R^2_{\gamma_i}(x_i,F(x_i))]
\leq \frac{8}{\alpha(k+1)}\Big[V(x_0,x^*)+\frac{C_4(x^*)+MC_5(x^*)}{\mathcal{N}b(\ln(\lambda-1))^b}\Big];$$
\item [(ii)]  If $N_k=\lceil \frac{(k+1)^{1.5}}{d}\rceil$, then
$$\min_{i\in\{0,\cdots,k\}}\mathbb{E}[R^2_{\gamma_i}(x_i,F(x_i))]
\leq\frac{8}{\alpha(k+1)}\Big[V(x_0,x^*)+\frac{C_4(x^*)+MC_5(x^*)}{3d}\Big].$$
\end{itemize}
\end{theorem}

\noindent{\bf Proof.} (i)
Taking expectation in \eqref{equ:4:18}, using $\mathbb{E}[\mathbb{E}[\cdot|\mathcal{F}_i]]=\mathbb{E}[\cdot]$ and summing it from $i=0$ to $i=k$, By Lemma \ref{lemma:4.6},
\begin{align*}
\frac{\alpha}{8}\sum^{k}_{i=0}\mathbb{E}[R^2_{\gamma_i}(x_i,F(x_i))]&\leq V(x_0,x^*)+C_4(x^*)\sum^{k}_{i=0}\frac{1}{N_i}+C_5(x^*)\sum^{k}_{i=0}\frac{\mathbb{E}[V(x_i,x^*)]}{N_i}\\
&\leq V(x_0,x^*)+C_4(x^*)\sum^{k}_{i=0}\frac{1}{N_i}+MC_5(x^*)\sum^{k}_{i=0}\frac{1}{N_i}.
\end{align*}
By dividing both sides of the above inequality by $k+1$, we obtain
\begin{align*}
\min_{i\in\{0,\cdots,k\}}\mathbb{E}[R^2_{\gamma_i}(x_i,F(x_i))]\leq& \frac{1}{k+1}\sum^{k}_{i=0}\mathbb{E}[R^2_{\gamma_i}(x_i,F(x_i))]\\
\leq&\frac{8}{\alpha(k+1)}\Big[V(x_0,x^*)+C_4(x^*)\sum^{k}_{i=0}\frac{1}{N_i}+MC_5(x^*)\sum^{k}_{i=0}\frac{1}{N_i}\Big].
\end{align*}
From \eqref{equ:4:14} and \eqref{equ:4:15}, we get the conclusions (i) and (ii) immediately. $\hfill\square$

For the unbounded case, the proof of the complexity is similar to Corollary \ref{cor:4.1}, so we omit its proof.

\begin{cor} \label{cor:4.2}
Let Assumptions \ref{ass:2.1} and \ref{ass:3.1}-\ref{ass:3.3}  hold. For any $x^*\in X^*$, $k\geq0$, set $N_k=\mathcal{N}\lceil (k+\lambda)(\ln (k+\lambda))^{1+b}\rceil$ or
$N_k=\lceil \frac{(k+1)^{1.5}}{d}\rceil$. Then the following statements are true:
\begin{itemize}
\item [(i)] In terms of $\min\limits_{i\in\{0,\cdots,k\}}\mathbb{E}[R^2_{\gamma_i}(x_i,F(x_i))]$, the iteration complexity of Algorithm 1 is $\mathcal{O}(\frac{1}{\varepsilon}$);

\item [(ii)] Both in the sense of almost surely and in the sense of expectation, the oracle complexity of Algorithm 1 is $\mathcal{O}(\frac{1}{\varepsilon^2})$.
\end{itemize}
\end{cor}

\begin{remark}
In Theorem \ref{the:4.2}, we obtain the convergence rate of the algorithm with the help of the Minty variational inequality, which is weaker than the strongly pseudomonotonicity used in \cite{KS}.
\end{remark}

\section{Numerical experiments}\label{sec:Numerical}
\noindent In this section, we give two numerical examples to illustrate the performance of Algorithm 1. In Section \ref{sec:Numerical-1}, we compare Algorithm 1 with the Variance-based extragradient algorithm in \cite{I2019} (shortly EGLS), the Mirror prox stochastic approximation algorithm in \cite{KS} (shortly by MPSA ) and the Variance-based subgradient extragradient algorithm in \cite{YZWL} (shortly by SEGVR). We also demonstrate the influence of different distance metrics on the numerical performance of Algorithm 1. In Section \ref{sec:Numerical-2}, we apply Algorithm 1 to solve the stochastic Nash game problem and give corresponding numerical analysis. All algorithms are coded in MATLAB R2019a and run the same computer with Windows 10 system, AMD Ryzen 5 3550H with Radeon Vega Mobile Gfx 2.10 GHz and 16.0 GB.

\subsection{Comparing Algorithm 1 with other algorithms}\label{sec:Numerical-1}
\noindent In this subsection, we compare Algorithm 1 with the algorithms EGLS, MPSA, SEGVR in terms of the CPU time and the residual value (denoted by VRF below) through averaging across 20 sample paths. For brevity, we use ``Iteration'' to refer to the maximum number of iterations, and denote by VRF:=$\|x_k-\mathbf{P}_X(x_k,\theta^{-1}\gamma_0\widehat {F}(x_k,\xi_k))\|$.  In our tests, all algorithms are terminated when the total number of iterations reaches the maximum iteration $K$.
Other parameters are chosen as follows:

Algorithm 1$: \gamma_0=0.99, \theta= 0.01, \alpha=2, \mathcal{Q}=5;$

EGLS$: \gamma_0=0.99, \theta= 0.01, \alpha=0.5;$

MPSA$: \gamma_0=0.99, \gamma_k=\frac{\gamma_0}{k};$

SEGVR$: \gamma_0=0.99, \theta= 0.01, \alpha=0.5.$

\begin{ex} \label{Example1} \cite{KS} (Stochastic Fractional Convex Quadratic Problem)
Consider the stochastic fractional convex quadratic problem of the form:
\begin{align*}
&\min_{x\in \mathbb{R}^n}\mathbb{E}\Big[\frac{\phi(x,\xi)}{\psi(x)}\Big]\\
& s.t. \, Ax\leq b,\; \left\|x\right\|_{1}\leq 1,
\end{align*}
where
\begin{align*}
&\phi(x,\xi):=0.5x^T\Big(0.025UU^T+0.025\frac{\left\|UU^T\right\|_F}{\left\|V(\xi)\right\|_F}V(\xi)\Big)x+0.5((c+\bar c(\xi))^Tx+4n)^2,\\
&\psi(x):=r^Tx+t+4n,
\end{align*}
$A\in \mathbb{R}^{m\times n}$, $b\in \mathbb{R}^{m\times 1}$, $m:=\lceil n/10\rceil$ and $\left\|\cdot\right\|_F$ denotes the Frobenius norm. Besides, $V(\xi)$ and $\bar c(\xi)$ are randomly generated from standard normals and uniform distributions, $U$ and $c$ are deterministic constants generated from standard normal distributions on [0,1], while $r$ and $t$ are generated from uniform distributions on [0,5].
\end{ex}

\begin{ex} \label{Example2}  (Fractional Convex Nonlinear Problem)
We consider a nonlinear variant of Example \ref{Example1} with same parameters and numerator but an exponential denominator $$\psi(x):=10^4\Big(e^{(8n+2)/2000}-e^{(r^Tx+t+4n)/2000}\Big).$$
\end{ex}

\begin{figure}[H]
\subfigure[$n$=10]{
\begin{minipage}[t]{0.5\textwidth}
		\centering
		\includegraphics[width=8cm]{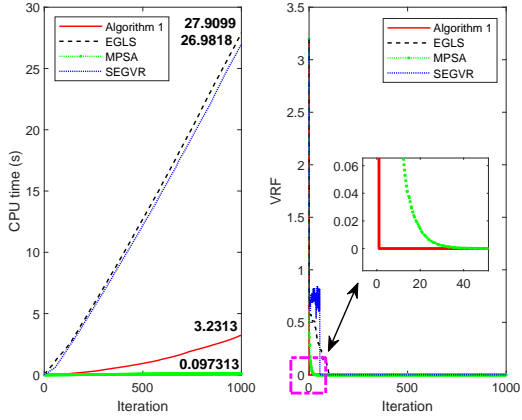}
		%\noindent{\bf (a)}\;{$n$=10}%\caption*{$n$=10}
	\end{minipage}
}
	\qquad
\subfigure[$n$=50]{
	\begin{minipage}[t]{0.5\textwidth}
		\centering
		\includegraphics[width=8cm]{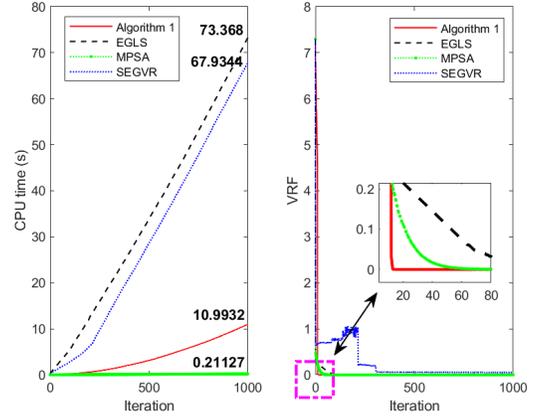}
		%\noindent{\bf (b)}\;{$n$=50}%\caption*{$n$=50}
	\end{minipage}
}
\caption{Results for Example \ref{Example1}  with $N_{k1}$.}
\label{fig1}
\end{figure}

\begin{figure}[H]	
\subfigure[$n$=10]{
\begin{minipage}[t]{0.5\textwidth}
		\centering
		\includegraphics[width=8cm]{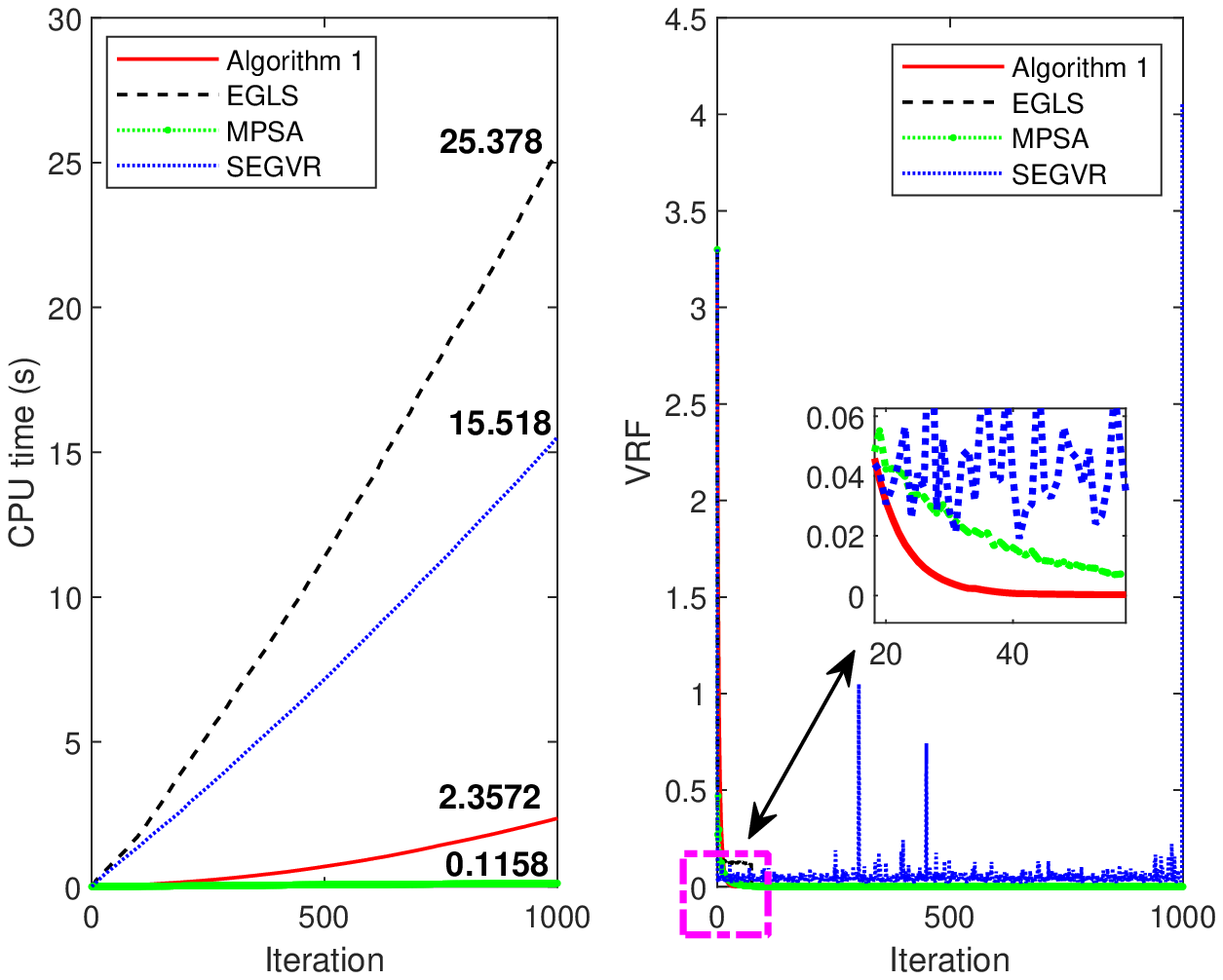}
		%\noindent{\bf (a)}\;{$n$=10}
	\end{minipage}
}
	\qquad
\subfigure[$n$=50]{
	\begin{minipage}[t]{0.5\textwidth}
		\centering
		\includegraphics[width=8cm]{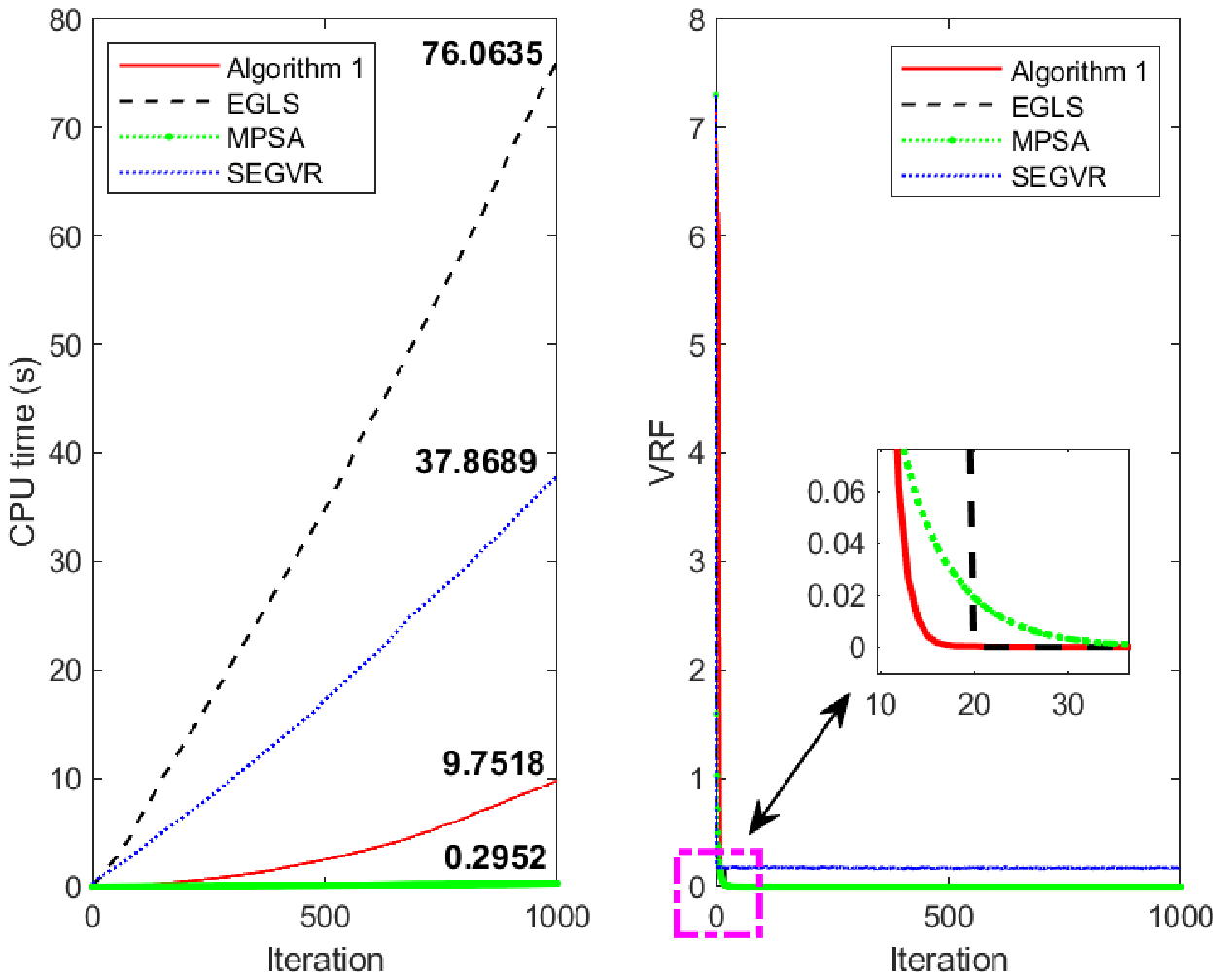}
		%\noindent{\bf (b)}\;{$n$=50}
	\end{minipage}
}
\caption{Results for Example \ref{Example2} with $N_{k1}$.}
\label{fig2}
\end{figure}

\begin{figure}[H]	
\subfigure[$n$=10]{
\begin{minipage}[t]{0.5\textwidth}
		\centering
		\includegraphics[width=8cm]{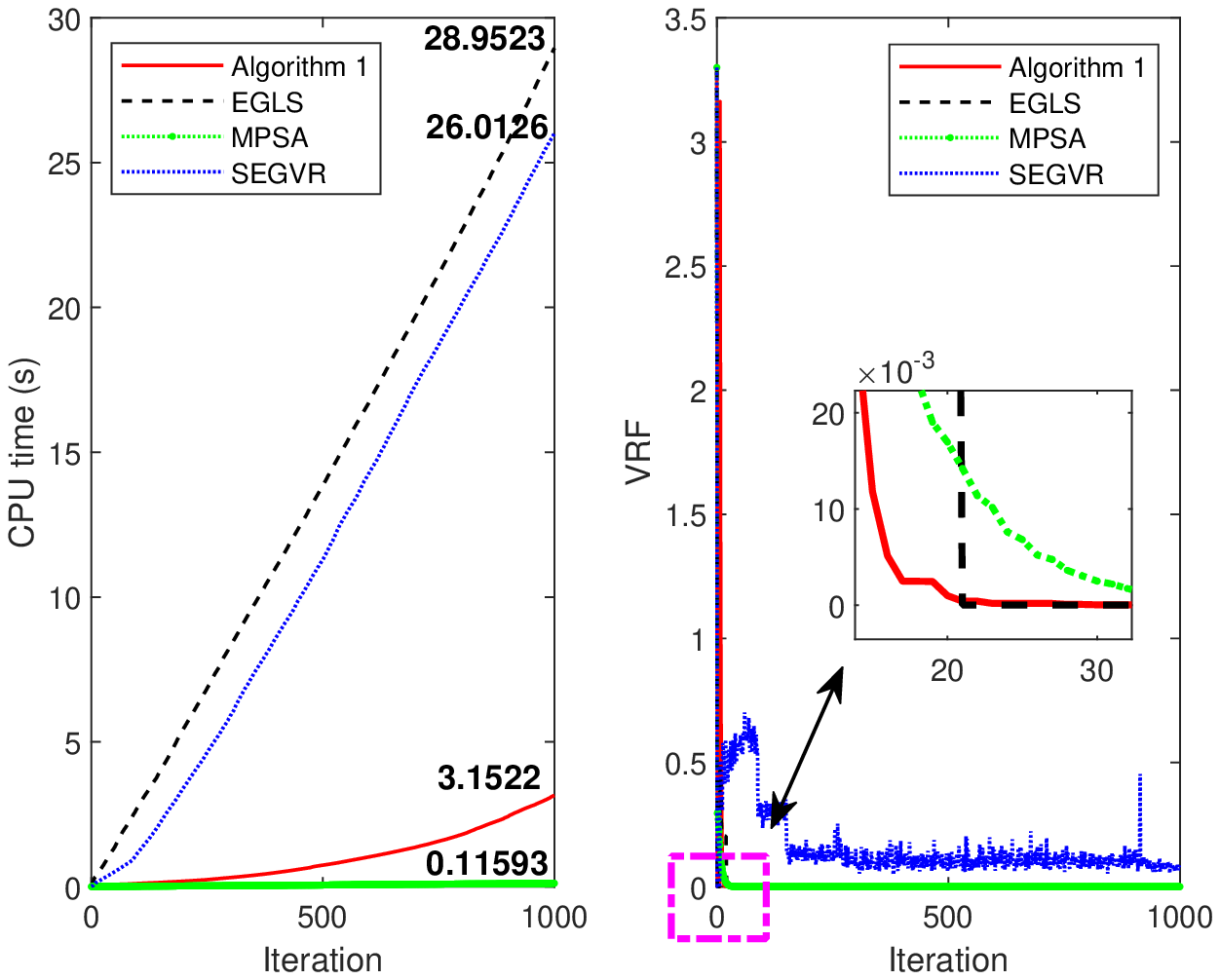}
		%\noindent{\bf (a)}\;{$n$=10}
	\end{minipage}
}
	\qquad
\subfigure[$n$=50]{
	\begin{minipage}[t]{0.5\textwidth}
		\centering
		\includegraphics[width=8cm]{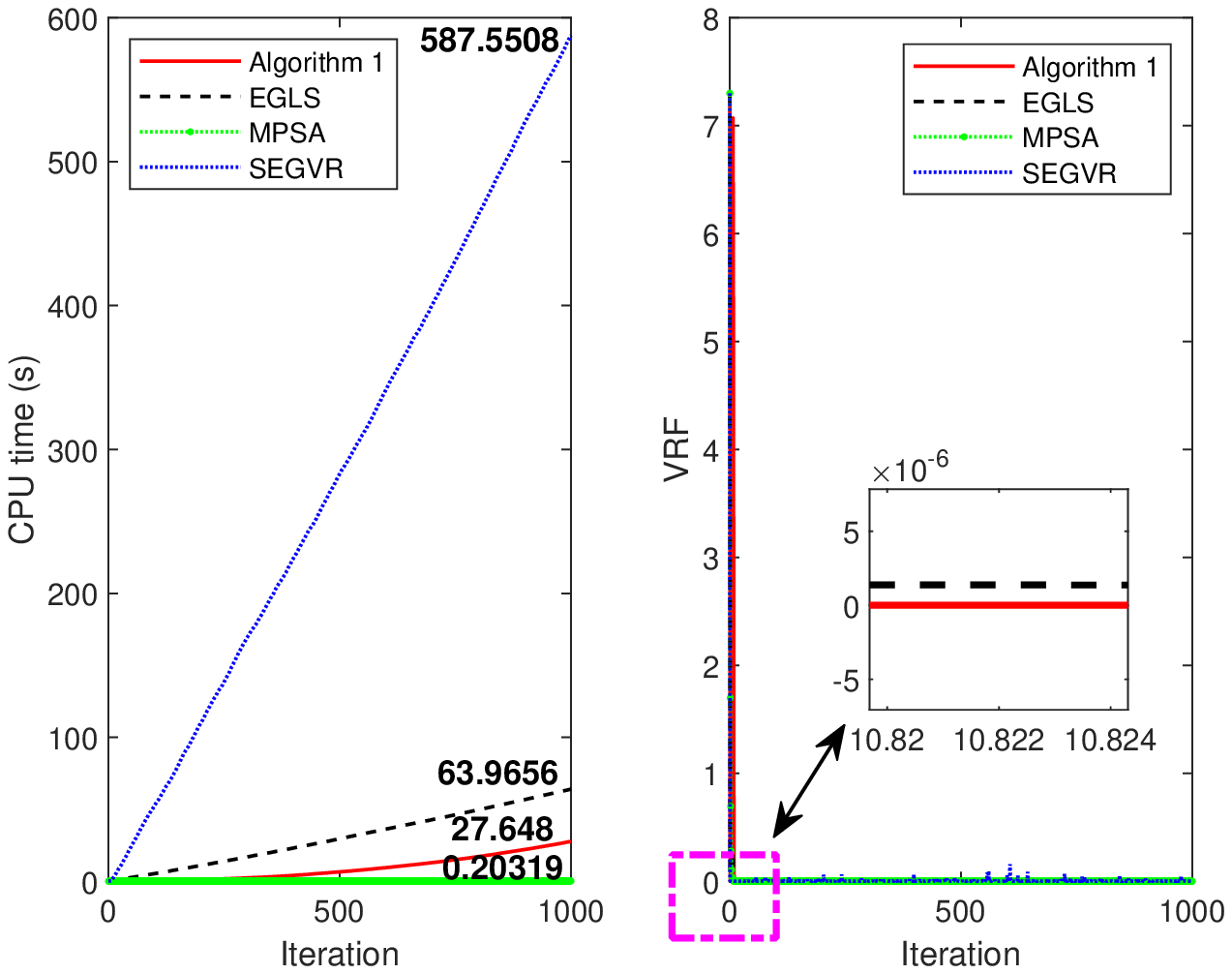}
		%\noindent{\bf (b)}\;{$n$=50}
	\end{minipage}
}
\caption{Results for Example \ref{Example1} with $N_{k2}$.}
\label{fig3}
\end{figure}

\begin{figure}[H]	
\subfigure[$n$=10]{
\begin{minipage}[t]{0.5\textwidth}
		\centering
		\includegraphics[width=8cm]{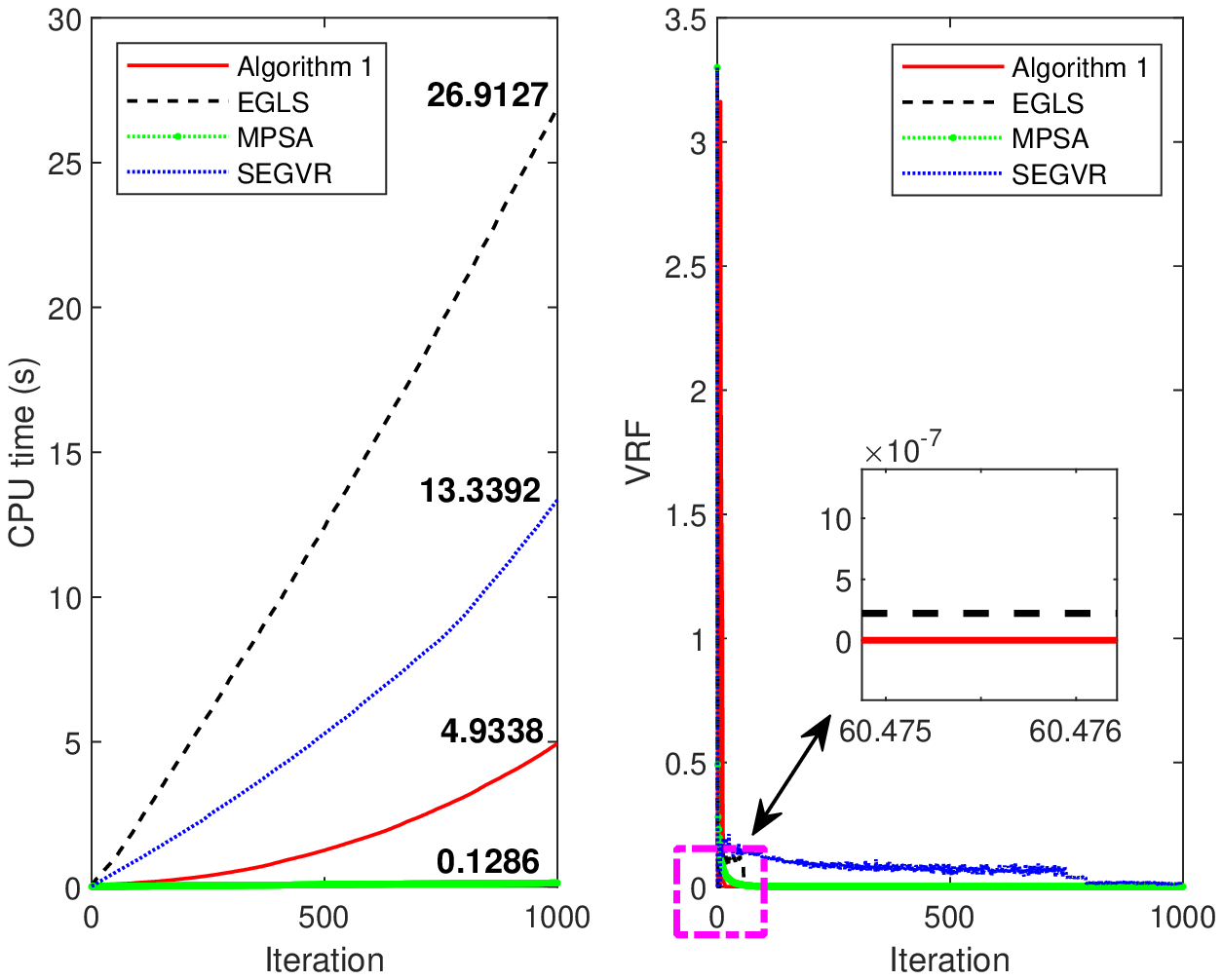}
		%\noindent{\bf (a)}\;{$n$=10}
	\end{minipage}
}
	\qquad
\subfigure[$n$=50]{
	\begin{minipage}[t]{0.5\textwidth}
		\centering
		\includegraphics[width=8cm]{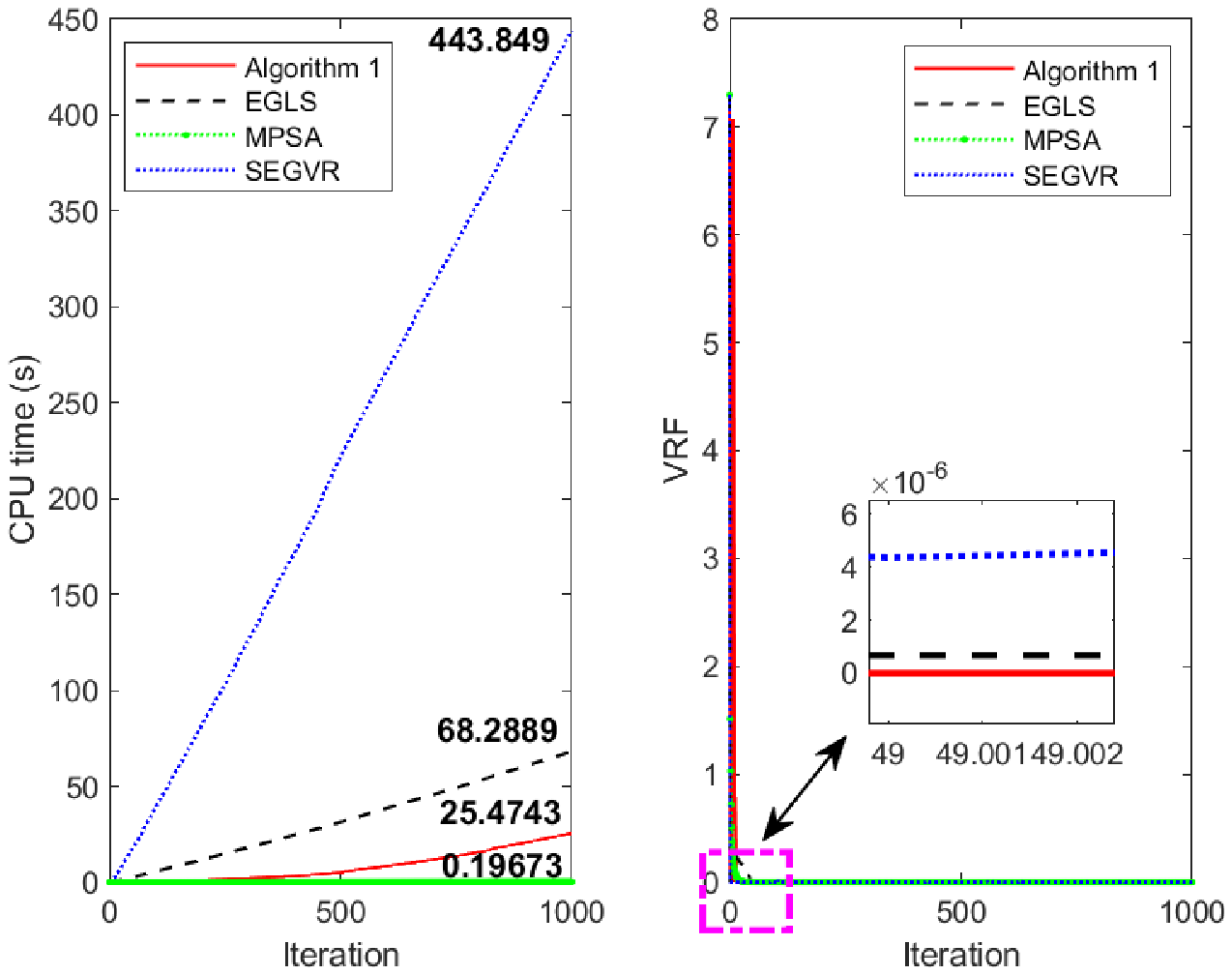}
		%\noindent{\bf (b)}\;{$n$=50}
	\end{minipage}
}
\caption{Results for Example \ref{Example2} with $N_{k2}$.}
\label{fig4}
\end{figure}

Firstly, we show the numerical comparison among Algorithms 1, EGLS, MPSA and SEGVR by means of the CPU time and the residual value with different samples $N_k$. We set $N_{k1}:=\lceil (k+1)^{0.8}\rceil$ and $ N_{k2}:=\frac{1}{10}\lceil (k+\lambda)(\ln (k+2.05))^{1.0001}\rceil$ and denote by the maximum iteration $K$=1000. The numerical results are shown in Figs.\ref{fig1}-\ref{fig4}.

%\noindent{\bf Fig.1} Results for Example \ref{Example1}  with $N_{k1}$.

%\noindent{\bf Fig.2} Results for Example \ref{Example2} with $N_{k1}$.

%\noindent{\bf Fig.3} Results for Example \ref{Example1} with $N_{k2}$.

%\noindent{\bf Figure.4} Results for Example \ref{Example2} with $N_{k2}$.\\

\renewcommand{\arraystretch}{1.5} %¿ØÖÆÐиß
%\begin{table}[tp] ±í¸ñÖö¥
\begin{table}[htb]
  \centering
  \fontsize{9}{10}\selectfont
  \begin{threeparttable}
  \caption{Algorithm $1$ with different distance metrics.}
  \label{table1}
    \begin{tabular}{lccccccc}
    \toprule
    \multirow{2}{*}{}&
    \multicolumn{2}{c}{$s_1(x)$}&\multicolumn{2}{c}{$s_2(x)$ }&\multicolumn{2}{c}{$s_3(x)$}\cr
    \cmidrule(lr){2-3} \cmidrule(lr){4-5} \cmidrule(lr){6-7}
    &CPU time (s) &VRF  &CPU time (s) &  VRF&   CPU time (s) &VRF\cr
    \midrule
    Ex.1($n$=10)\cr
    K=500    &0.863   &2.478e$-$11   &0.652    &2.000e$-$03    &1.008    &5.175e$-$05\cr
    K=1000   &3.413   &6.737e$-$12   &2.251    &1.000e$-$03    &5.149    &3.381e$-$05\cr
    Ex.1($n$=50)\cr
    K=500    &3.493   &1.885e$-$10   &2.450    &3.960e$-$01    &3.850    &8.174e$-$04\cr
    K=1000   &12.097  &9.655e$-$13   &7.523    &3.200e$-$02    &13.912   &7.806e$-$04\cr

    Ex.2($n$=10)\cr
    K=500    &1.039   &1.484e$-$14   &0.801    &3.300e$-$02    &0.892    &5.324e$-$04\cr
    K=1000   &3.573   &1.952e$-$15   &2.658    &3.200e$-$02    &2.982    &3.523e$-$04\cr
    Ex.2($n$=50)\cr
    K=500    &2.875   &1.687e$-$11   &2.388    &3.050e$-$01    &2.365    &5.400e$-$03\cr
    K=1000   &9.722   &6.793e$-$12   &8.727    &2.620e$-$01    &9.010    &2.500e$-$03\cr
    \bottomrule
    \end{tabular}
    \end{threeparttable}
\end{table}

In Figs.\ref{fig1}-\ref{fig4}, we show the performance of all algorithms for Examples \ref{Example1} and \ref{Example2}. It can be observed that the values of VRF solved by Algorithm 1 vary with the same decreasing trend as Algorithms EGLS and MPSA, but our algorithm tends to 0 more quickly than Algorithms EGLS, MPSA and SEGVR. Since Algorithm MPSA only adopts a single sample strategy, it spends a little CPU time but has the high complexity. However, Algorithm 1 saves at least 80$\%$ CPU time in comparison with Algorithm EGLS and enjoys lower complexity.
As for the sample size, it is easy to see that for either $N_{k1}$ or $N_{k2}$, Algorithm 1 always performs best among all compared algorithms in terms of the values of VRF. Moreover, compared with the setting of $N_{k2}$, Algorithm 1 spends less CPU time when taking the sample size as $N_{k1}$.
To sum up, Algorithm 1 is competitive with Algorithms EGLS, MPSA and SEGVR.

Secondly, as mentioned before, one salient advantage of the Bregman distance is that they can change with different distance functions, which may deal with more general problems in practice. Hence, we consider three different distance functions given by \cite{NA} to illustrate this fact. For convenience, we show the CPU time and the VRF values of Algorithm 1 with respect to Examples \ref{Example1} and \ref{Example2}, where the distance generating functions are defined by $s_1(x):=\sum^{n}_{i=1}(x_i+\sigma)\log (x_i+\sigma)$, $s_2(x):=\frac{1}{2}\sum^{n}_{i=1}\|x_i\|^2$ and $s_3(x):=\log (n)\sum^{n}_{i=1}x_i^{\big(1+\frac{1}{\log(n)}\big)}$, respectively. The specific numerical results are shown in Table \ref{table1}.

As revealed in Table \ref{table1}, there are obvious distinctions on the numerical performance of Algorithm 1 when using different distance functions. In light of the CPU time, Algorithm 1 has slightly larger advantage of taking $s_2(x)$, in comparison with $s_1(x)$ and $s_3(x)$. In light of the VRF values, the VRF values solved by Algorithm 1 with $s_1(x)$ performs best, followed by $s_3(x)$ and $s_2(x)$. Therefore, Algorithm 1 with $s_1(x)$ is better than that with $s_2(x)$ and $s_3(x)$ overall.

\subsection{Application of Algorithm 1 to stochastic Nash game}\label{sec:Numerical-2}
\noindent Consider a networked Nash-Cournot game with uncertain data in \cite{KS,KNS,KS1,YSS}. Assume that there are $I$ firms (players), denoted by $\mathcal{I}:=\{1,\cdots,I\}$, that compete over a network of $J$ spatially distributed markets (nodes), denoted by $\mathcal{J}:=\{1,\cdots,J\}$. Firm $i$ needs to determine a continuous-valued nonnegative quantities, denoted by $x_i:=(x^1_i, \cdots, x^J_i)^T\in\mathbb{R}^J$, to be delivered to the markets, where $x^j_i$ represents the sales of firm $i$ at market $j$. Let $\bar {x}^j:=\sum_{i\in \mathcal{I}} x^j_i$ denote the aggregate sales at market $j$. By the Cournot structure, we assume that the price $p_j$ sold in market $j$ is determined by a linear inverse demand function corrupted by $p_j(\bar x^j,\xi):=a_j(\xi)-b_j\bar x^j$, where $a_j(\xi)$ is a uniform random variable drawn from $[l_{a_j(\xi)},u_{a_j(\xi)}]$ and the nonnegative parameter $b_j$ represents the slope of the inverse demand function. Moreover, firm $i$ has a random linear linear production cost function $\mathbb{C}_i(x_i,\xi):=c_i(\xi)\sum_{j\in\mathcal{J}}x^j_i$, where $c_i(\xi)$ is a uniform random variable drawn from $[l_{c_j(\xi)},u_{c_j(\xi)}]$. Suppose that firm $i$ has finite production capacity, i.e.,
$$X_i:=\{x_i\in \mathbb{R}^J|x_i\geq 0,x^j_i\leq cap^j_i,\forall j\in\mathcal{J}\}.$$
It follows that firm $i$'s optimization problem is given by
\begin{align*}
\min_{x_i\in X_i}\mathbb{E}[\theta_i(x,\xi)]:&=\mathbb{E}\Big[\mathbb{C}_i(x_i,\xi)-\sum_{j\in J}p_j(\bar x^j,\xi)x^j_i\Big]\\
&=\mathbb{E}\Big[x^T_idiag(b)x_i+\big(c_i(\xi)e+\sum_{s\neq i, s\in \mathcal{I}}b\circ x_s-a(\xi)\big)^Tx_i\Big],
\end{align*}
where $\circ$ denotes the Hadamard product, $e$ denotes the vector whose elements are all one, $a(\xi):=(a_1(\xi),\cdots,a_J(\xi))^T$ and $b:=(b_1,\cdots,b_J)^T$. It is easy to see that this problem is convex. Therefore, by the optimization theory, the stochastic Nash-Cournot game can be compactly formulated as a SVI with $F(x^*):=(F_1(x^*),\cdots,F_I(x^*)), F_i(x^*)=\mathbb{E}[\partial_{x_i}\theta_i(x,\xi)]$, and $X := \Pi_{i\in \mathcal{I}}X_i$.

In our tests, let $\mathcal{I}=10,20,30$ and $\mathcal{J}=10$. For any $i\in \mathcal{I}$ and $j\in \mathcal{J}$, let $[l_{a_j(\xi)}$, $u _{a_j(\xi)}]=[30,60]$, $[l_{c_j(\xi)},u_{c_j(\xi)}]=[2,6]$, $cap^j_i=2$ and $b_j$ be a deterministic constant generated once from uniform distribution on $[0,2]$. Set $N_K=2\lceil (K+1)^{0.8}\rceil$ with Relative error:$=\frac{\|x^K-x^*\|}{\|x^*\|}$. The numerical results are shown in Table \ref{table2}.

\renewcommand{\arraystretch}{1.5} %¿ØÖÆÐиß
%\begin{table}[tp] ±í¸ñÖö¥
\begin{table}[htb]
  \centering
  \fontsize{9}{9.5}\selectfont
  \begin{threeparttable}
  \caption{The numerical results of the stochastic Nash game}
  \label{table2}
    \begin{tabular}{llll}
    \toprule
    The number of the firm&Iteration (K)&Relative error&CPU time (s)\cr
    \midrule
   I=10 &100	    &1.342e$-$01	&9.210e$-$02\cr
        &500	    &4.070e$-$02	&9.443e$-$01\cr
        &1000       &5.000e$-$03	&2.961e+00\cr
        &2000       &2.500e$-$03	&1.087e+01\cr
        &5000       &9.793e$-$04    &5.479e+01\cr

   I=20 &100	    &1.072e$-$01	&9.160e$-$02\cr
        &500	    &3.160e$-$02	&1.040e+00\cr
        &1000	    &4.200e$-$03	&3.248e+00\cr
        &2000       &2.400e$-$03	&1.005e+01\cr
        &5000       &8.616e$-$04    &5.632e+01\cr

   I=30 &100	    &1.041e$-$01	&8.680e$-$02\cr
        &500	    &2.910e$-$02	&1.067e+00\cr
        &1000       &1.000e$-$02	&3.317e+00\cr
        &2000       &3.600e$-$03	&1.282e+01\cr
        &5000       &8.360e$-$04    &6.312e+01\cr
    \bottomrule
    \end{tabular}
    \end{threeparttable}
\end{table}

From Table \ref{table2}, it is easy to see that increasing the iteration $K$ and the number of firm $I$ may slow down the values of relative error, while increase more CPU time. Besides, the proposed algorithm may obtain approximation solutions which are very close to the ground truth solutions when $K=5000$.

\section{Conclusion}\label{sec:conclusion}
\setcounter{equation}{0}

\noindent In this paper, we propose a new stochastic Bregman extragradient algorithm with line search for solving the SVI. We prove almost sure convergence of the algorithm in terms of the natural residual function. Under the Minty variational inequality condition, we obtain the convergence rate with respect to the gap function and the natural residual function, respectively, as well as the iteration complexity and the oracle complexity. Finally, several numerical experiments illustrate the efficiency and competitiveness of the proposed algorithm. As for the future work, we will further study the Bregman single projection algorithm so as to speed up the convergence. Extending the proposed method to stochastic mixed variational inequalities is also one of our future goals.

{\footnotesize{
}}
\end{document}